\documentclass[12pt, reqno]{amsart}
\usepackage{amsfonts}
\usepackage{amsmath}
\usepackage[latin1]{inputenc}
\usepackage[english]{babel}
\usepackage{amsthm}
\usepackage{graphicx}
\usepackage{color}

\usepackage{amssymb}
\usepackage{amsopn}
\usepackage{enumerate}
\usepackage{calrsfs}
\usepackage{comment}

\usepackage{pifont,bbding}
\usepackage{comment}
\usepackage{mathrsfs}
\usepackage{cancel}

\newcommand{\N}{\mathbb{N}}
\newcommand{\R}{\mathbb{R}}
\newcommand{\Rn}{\mathbb{R}^n}

\newcommand{\F}{\mathcal{F}}
\newcommand{\HH}{\mathcal{H}}

\newcommand{\diva}{\mathrm{div}\!_\alpha}
\renewcommand{\P}{P_\alpha}
\newcommand{\PH}{P_H}
\newcommand{\PX}{P_X}
\renewcommand{\L}{\mathcal L}
\newcommand{\I}{\mathcal I}

\newcommand{\Q}{Q}

\renewcommand{\S}{\mathcal S_x}

\theoremstyle{plain}
\newtheorem{thm}{Theorem}[section]
\newtheorem{prop}[thm]{Proposition}

\theoremstyle{definition}

\theoremstyle{remark}
\newtheorem{rem}[thm]{Remark}

\numberwithin{equation}{section}

\title[Isoperimetric problem]
{Isoperimetric problem in $H$-type groups \\ and Grushin spaces}

\long\def\MSC#1\EndMSC{\def\arg{#1}\ifx\arg\empty\relax\else
     {\par\narrower\noindent
     {\small\it 2010 Mathematics Subject Classification.} \small #1\par}\fi}

\long\def\KEY#1\EndKEY{\def\arg{#1}\ifx\arg\empty\relax\else
     {\par\narrower\noindent
     {\small\it Keywords and Phrases.} \small #1\par}\fi}

\author[V.~Franceschi]{Valentina Franceschi}
\email{valentina.franceschi@math.unipd.it}

\author[R.~Monti]{Roberto Monti}
\email{monti@math.unipd.it}

\address[Franceschi and Monti]
{Universit\`a di Padova, Dipartimento di Matematica,
via Trieste 63, 35121 Padova, Italy}

\subjclass[2010]{53C17,49Q20}
\keywords{Isoperimetric Problem, $H$-type groups, Grushin spaces.}

\thanks{This work was partially supported by the Fondazione CaRiPaRo Project
``Nonlinear Partial Differential Equations: models, analysis, and
control-theoretic problems'', Padova.}

\date{\today}

\begin{document}

\begin{abstract}
We study the isoperimetric problem in $H$-type groups and 
Grushin spaces, emphasizing a relation between them. 
We prove existence, symmetry and regularity properties of
isoperimetric sets, under a symmetry assumption that depends
on the dimension.
\end{abstract}

\maketitle

\section{Introduction}

Let $M$ be a manifold,   $V$ be a volume, and $P$ a perimeter
measure on $M$. For a   regular set   $E\subset M$, $P(E)$
is the  area of the boundary $\partial
E$. The isoperimetric    
problem relative to $V$ and $P$  
consists in studying
existence,   
symmetries, regularity and, if possible, classifying the minimizers
of the   problem
\begin{equation}\label{Isop}
    \min \big\{ P(E) : E \in \mathcal A  \text{ such that  } V(E)= v
\big\},
\end{equation}
for a given volume $v>0$ and for a given family of admissible sets $\mathcal
A  $. Minimizers of \eqref{Isop} are called isoperimetric
sets.

In space forms (Euclidean space, sphere and hyperbolic space) with their natural
volume and perimeter, isoperimetric sets 
are precisely metric balls. In $\Rn$ with volume 
 $\mathrm{e}^{-|x|^2} \mathcal L^n $
and perimeter 
$\mathrm{e}^{-|x|^2} \mathcal H^{n-1} $,
isoperimetric sets are  half-spaces. This is the Gaussian isoperimetric
problem, the model of the current research direction on
isoperimetric problems with density. A different way to wheight  
perimeter is by a surface tension, i.e., by  the support function $\tau:\mathbb
S^{n-1}\to[0,\infty)$
of a convex body $K\subset \Rn$ with $0\in\mathrm{int}(K)$, $\tau(\nu) =
\sup_{x\in K}\langle x,\nu\rangle$.  Namely, one can
consider
\[
     P(E) = \int_{\partial E}  \tau(\nu_E) d\mathcal H^{n-1},\quad \nu_E \text{
outer normal to $\partial E$.}
\]
The isoperimetric problem for this perimeter and with $V=\mathcal L^n$  is known
as Wulff problem and isoperimetric sets are translates and dilates of the
set $K$.

In a different approach, the perimeter of a Lebesgue measurable set 
$E\subset\Rn$   is defined via a system  
$X=\{ X_1,\ldots,X_h\}$, $h\geq
2$,  of self-adjoint vector fields in $\Rn$, $X_j=-X_j^*$,
\begin{equation}\label{pipp}
 \PX(E) =   \sup \Big\{ \int_{E} \sum_{i=1}^h X_i \varphi_i(x) \,
dx \, :\, \varphi\in C_c^1(\Rn;\R^{h}),\, 
\displaystyle  \max_{x\in \Rn } |\varphi(x)|  \leq 1
 \Big\}.
\end{equation}
This definition is introduced and  studied systematically in \cite{GN}. The
perimeter $\PX$ is known as $X$-perimeter (horizontal, sub-elliptic, or
sub-Riemannian perimeter). One
important example
is the Heisenberg  perimeter, that is subject of intensive research in
connection with Pansu's conjecture
on the shape of  
isoperimetric sets (see  \cite{MR}, \cite{RR1},
\cite{RR2}, \cite{R}, \cite{M3}) and  in connection with the
regularity problem of minimal surfaces.

In this paper, we study perimeters that are related to   the Heisenberg
perimeter.
Namely, we study the isoperimetric problem  
in $H$-type  groups and  in Grushin spaces.

\medskip

1) {\emph{ $H$-type groups.}} Let $\mathfrak h = \mathfrak h_1\oplus\mathfrak
h_2$
be a stratified nilpotent real Lie algebra of dimension $n\geq 3$ and step $2$.
Thus we have $\mathfrak h_2 =[\mathfrak h_1,\mathfrak h_1]$. We fix on
$\mathfrak h$
a
scalar product $\langle\cdot,\cdot\rangle$ that makes $\mathfrak h_1$ and
$\mathfrak h_2$ orthogonal. 
The Kaplan mapping is the mapping $J:\mathfrak h_2\to \mathrm{End}(\mathfrak
h_1)$ defined via the identity
\begin{equation}\label{Kaplan}
 \langle J_Y (X), X'\rangle = \langle Y, [X,X']\rangle,
\end{equation}
holding for all $X,X'\in \mathfrak h_1$ and $Y\in\mathfrak h_2$. The algebra
$\mathfrak
h$ is
called an $H$-type algebra if for all $X,X'\in\mathfrak h_1$ and $Y\in\mathfrak
h_2$ there holds
\begin{equation}\label{H-type}
 \langle J_Y(X),J_Y(X')\rangle = |Y|^2 \langle X,X'\rangle,
\end{equation}
where $|Y|= \langle Y,Y\rangle^{1/2}$.
We can identify
$\mathfrak h$ with $\R^n = \R^h\times\R^k$, $\mathfrak h_1$ with
$\R^h\times\{0\}$,
and $\mathfrak h_2$ with $\{0\}\times \R^k$, where
$h\geq 2$ and $k\geq 1$ are integers.
In fact, $h$ is an even integer.
We can also assume that
$\langle\cdot,\cdot\rangle$ is the standard scalar product of $\R^n$.
Using exponential coordinates, 
the connected and simply connected
Lie group of $\mathfrak h$ can be identified with $\Rn$.
Denoting points of $\Rn$
as $(x,y) \in\R^n=\R^h\times\R^k$, 
the Lie group product $\cdot :\Rn\times\Rn\to\Rn$
  is of the form $(x,y) \cdot (x',y') = (x+x',y+y'+Q(x,x'))$, 
where $Q:\R^h\times\R^h\to\R^k$ is a bilinear skew-symmetric mapping. 
Let $Q_{ij}^\ell\in\R$ be the numbers
\[
   Q^\ell_{ij}  = \langle Q
(\mathrm{e}_{i},\mathrm{e}_{j}),\mathrm{e}_\ell\rangle, \quad
i,j=1,\ldots,h,\,\,
\ell=1,\ldots,k,
\]
where $\mathrm{e}_{i},\mathrm{e}_{j}\in\R^h$ and $\mathrm{e}_\ell\in\R^k$ are
the standard coordinate  versors.
An orthonormal basis of the Lie algebra of left-invariant vector fields of the
$H$-type group $(\Rn,\cdot)$ is
given
by
\begin{equation}\label{HTYP}
\begin{split}
    X_i & = \frac{\partial}{\partial x_i} -\sum_{\ell=1}^k 
 \sum_{j=1}^h Q^\ell_{ij} x_j \frac{\partial}{\partial y_\ell},\quad
i=1,\ldots,h,
\\
   Y_j & = \frac{\partial}{\partial y_j},\quad j=1,\ldots,k.
\end{split}
\end{equation}
We denote by $\PH(E)=\PX(E)$ the perimeter of a set $E\subset\Rn$ defined as in
\eqref{pipp}, relatively to the system of vector fields $X=
\{X_1,\ldots,X_h\}$. 
The vector fields $Y_1,\ldots, Y_k$ are not considered.

\medskip

2) \emph {Grushin spaces.}
Let $\R^n = \R^h\times\R^k$, where   $h,k\geq 1$ are
integers and  $n = h+k$. 
For a given real  number $\alpha >0$, let us define the vector fields
in $\R^n$
\begin{equation}\label{GTYP}
\begin{array}{l}
\displaystyle X_i =\frac{\partial}{\partial {x_i} },
\quad i=1,\ldots,h,
\\
\displaystyle Y_j=|x|^{\alpha}\frac{\partial}{\partial {y_j} },
\quad j=1,\ldots,k,
\end{array}
\end{equation}
where $|x|$
is the standard norm of $x$.
We denote by $\P(E)=\PX(E)$ the perimeter of a set $E\subset\Rn$ defined as in
\eqref{pipp} relatively to the system of vector fields $X= \{X_1,\ldots,X_h,
Y_1,\ldots,Y_k\}$.
We call $\P(E)$ the $\alpha$-perimeter of $E$.

\medskip

We study the isoperimetric problem in the class of $x$-spherically symmetric
sets
in $H$-type groups and Grushin spaces. These two problems are related to
each other. We say that a set $E\subset \R^h\times\R^k$ is $x$-spherically
symmetric if there exists
a set $F\subset \R^+\times \R^k$, called  generating set of $E$, such that
\[
         E = \big\{(x,y)\in \R^n : (|x|,y) \in F\big\}.
\]
We denote by $\S$ the class of $\L^n$-measurable,
 $x$-spherically symmetric sets.

Starting from the $x$-spherical symmetry, 
we can prove that the class of sets involved in
the minimization \eqref{Isop} can be restricted to a smaller
class of  sets with more symmetries (see Section \ref{TRE}). 
Using this additional symmetry, we can implement the concentration-compactness
argument in order to have the existence of isoperimetric sets.
In Carnot groups, the existence is already known, see \cite{LR}.
In Grushin spaces, the existence is less clear because 
$x$-translations do not preserve $\alpha$-perimeter.

In fact, we
have existence of isoperimetric sets that are  $x$- and $y$-Schwartz symmetric,
i.e., of the form
\begin{equation}\label{SISSO}
E = \{(x,y)\in\Rn : |y|< f(|x|)\},
\end{equation} 
for some
function $f:(0,r_0)\to\R^+$,
 $r_0>0$, which is called 
the profile function of $E$. The profile function has the necessary
regularity to  solve  a  second order ordinary differential equation
expressing the fact that the boundary of $E$ has a certain ``mean curvature''
that is constant. 
This differential equation can be partially integrated and, 
for the profile function of a minimizer, 
it can be expressed in the following equivalent way:
\begin{equation}\label{ECCO}
   \frac{f'(r)}{\sqrt{r^{2\alpha} +f'(r)^2}} = \frac{k-1}{r^{h-1}}\int_0^r
\frac{s^{2\alpha+h-1}}{f(s)\sqrt{s^{2\alpha} +f'(s)^2} }ds -\frac{\kappa}{h}
r,\quad \text{for } r\in (0,r_0),
\end{equation}
where $h,k$ are the dimensional parameters, $\alpha>0$ is the real parameter in
the Grushin vector fields \eqref{GTYP}
(in $H$-type groups we have $\alpha=1$), and
$\kappa>0$
is a real parameter (the ``mean curvature'') related to perimeter and volume.

\medskip

In $H$-type groups, the Haar measure is the Lebesgue measure.
Moreover, Lebesgue measure and $H$-perimeter are homogeneous
with respect to the  anisotropic dilations
\[
(x,y)\mapsto \delta_\lambda(x,y)=(\lambda x,\lambda^{2}y),\quad \lambda>0.
\]
In fact, for any  measurable set $E\subset\R^{n}$
and   for all $\lambda>0$ we have $\L^n (\delta_\lambda(E))=\lambda^Q \L^n(E)$ and $\PH(\delta_\lambda(E))=\lambda^{Q-1}\PH(E)$, where the number $Q=h+2k $ is the 
homogenous dimension of the  group.
Then, the isoperimetric ratio
\[
\I_H(E) =\frac{\PH(E)^{Q}}{\L^n(E)^{Q-1}} 
\]
is homogeneous of degree $0$ and the isoperimetric problem
\eqref{Isop} can be formulated in scale invariant form.
In the following, by a vertical translation we mean a mapping of the form
$(x,y)\mapsto (x,y+y_0)$
for some $y_0\in\R^k$.

\begin{thm} \label{THM1} In any $H$-type group, the isoperimetric problem
\begin{equation}\label{IsopH}
    \min \big\{ \I_H (E) : E\in \S \textrm{ with }0<\L^n(E)<\infty \big\}
\end{equation}
has solutions and, up to a vertical translation and a null set,
any isoperimetric set is of the form \eqref{SISSO} for a function 
$f\in C([0,r_0])\cap C^1([0,r_0))\cap C^\infty(0,r_0)$, with $0< r_0 <\infty$,
satisfying $f(r_0)=0$, $f'\leq 0$ on $(0,r_0)$, and 
solving equation \eqref{ECCO} with $\alpha=1$ and $\kappa =
\frac{Q\PH(E)}{(Q-1) \L^n(E)}$.
\end{thm}

\noindent Isoperimetric sets are, in fact, $C^\infty$-smooth sets away
from  $y=0$. Removing the assumption of $x$-spherical symmetry is
a difficult problem that is open even in the basic example of the 
$3$-dimensional Heisenberg group.

\medskip
For the special dimension $h=1$, we are able to prove the $x$-symmetry
of isoperimetric sets for $\alpha$-perimeter. Lebesgue measure and
$\alpha$-perimeter are homogeneous
with respect to the group of anisotropic dilations
\[
(x,y)\mapsto \delta_\lambda(x,y)=(\lambda x,\lambda^{1+\alpha}y),
\quad \lambda>0.
\] 
In fact, for any  measurable set $E\subset\R^{n}$
and   for all $\lambda>0$ we have $\L^n (\delta_\lambda(E))=\lambda^d \L^n(E)$ and
$\P (\delta_\lambda(E))=\lambda^{d-1}\PH(E)$, 
where   $d=h+k(1+\alpha)$. Then, the isoperimetric ratio
\[
\I_\alpha(E) =\frac{\P (E)^{d}}{\L^n(E)^{d-1}} 
\]
is homogeneous of degree $0$.

\begin{thm} \label{THM2} Let $\alpha>0$, $ h =1$,  $k\geq 1$ and $n=1+k$.
The isoperimetric problem
\begin{equation}\label{IsopH1}
    \min \big\{ \I_\alpha (E) : E\subset \Rn \text{ $\L^n$-measurable with }
     0<\L^n(E)<\infty \big\}
\end{equation}
has solutions and, up to a vertical translation and a null set,
any isoperimetric set is of the form \eqref{SISSO} for a function 
$f\in C([0,r_0])\cap C^1([0,r_0))\cap C^\infty(0,r_0)$, with $0< r_0<\infty$,
satisfying $f(r_0)=0$, $f'\leq0$ on $(0,r_0)$, and 
solving equation \eqref{ECCO} with $h=1$ and $\kappa =
\frac{d\P (E)}{(d-1) \L^n(E)}$.
\end{thm}

 \noindent In particular, for $h=1$ isoperimetric sets are
$x$-symmetric. When $h\geq 2$ we need to assume the $x$-spherical symmetry.

\begin{thm} \label{THM3} Let $\alpha>0$, $ h \geq 2$,  $k\geq 1$ and $n=h+k$.
The isoperimetric problem
\begin{equation}\label{IsopH2}
    \min \big\{ \I_\alpha (E) : E\in\S   \text{ with  }
     0<\L^n(E)<\infty \big\}
\end{equation}
has solutions and, up to a vertical translation and a null set,
any isoperimetric set is of the form \eqref{SISSO} for a function 
$f\in C([0,r_0])\cap C^1([0,r_0))\cap C^\infty(0,r_0)$, with $0<r_0<\infty$,
satisfying $f(r_0)=0$, $f'\leq 0$ on $(0,r_0)$, and 
solving equation \eqref{ECCO} with $\kappa =
\frac{d\P (E)}{(d-1) \L^n(E)}$.
\end{thm}

\noindent In the special case $k=1$,   equation \eqref{ECCO}
can be integrated and we have an explicit formula for isoperimetric sets.
Namely, with the normalization $\kappa = h$ -- that implies  $r_0=1$, -- the
profile function solving  \eqref{ECCO} gives the isoperimetric set  
\begin{equation} \label{EXPO}
    E = \Big\{ (x,y) \in \Rn : |y|<\int _{\arcsin |x|}^{\pi/2} \sin^{\alpha+1}(s) \,
ds\Big\}.
\end{equation}
This formula  generalizes to dimensions $h\geq 2$ the results of \cite{MontiMorbidelli}.
When $k=1$ and $\alpha=1$, the profile function satisfying the final
condition $f(1)=0$ is $f(r) =\frac 12 \big(  \arccos(r) +r\sqrt{1-r^2}\big)$,
$r\in[0,1]$.
This is the profile function of the Pansu's ball in the Heisenberg group.

In Section \ref{DUE}, we prove various representation formulas for the
perimeter of smooth and symmetric sets. In particular, we show that
for $x$-spherically symmetric sets we have the identity $\PH(E) = \P(E)$ with
$\alpha=1$. This makes Theorem \ref{THM1} a special case of Theorem \ref{THM3}.

In Section \ref{TRE}, we prove the rearrangement theorems.
We show that when $h=1$ the isoperimetric problem with no symmetry assumption can be reduced to
$x$-symmetric sets. When $h\geq 2$, we show that the $x$-spherical symmetry can
be improved to the $x$-Schwartz symmetry. We also study perimeter under 
$y$-Schwartz rearrangement. The equality case in this rearrangement
does not imply that, before rearrangement, the set is already $y$-Schwartz
symmetric because the centers of the $x$-balls may vary. However,
for isoperimetric sets 
the centers are constant, see Proposition \ref{final}. To prove this, 
we use the regularity of the profile function (see Section \ref{QUATTRO}).

The existence of isoperimetric sets is established in Section
\ref{sec:existence} by the concentration-compactness method. Here, we 
borrow some ideas from \cite{FuscoMaggiPratelli}
and we also 
use the isoperimetric inequalities (with nonsharp
constants) obtained in \cite{GN}, \cite{FGW}, and \cite{FGuW}.

Finally, in Section \ref{QUATTRO} we deduce
the differential equation for the  profile function,  
we use minimality to derive its equivalent  version \eqref{ECCO},
and we establish some elementary properties of solutions.

\section{Representation and reduction formulas}\label{DUE}

In this section, we derive  some  formulas for the representation
of  $H$- and $\alpha$-perimeter of smooth sets
and of sets with symmetry.
For any open set $A\subset\Rn$ and $m\in\N$, let us define the family of test
functions
\[
 \F_m(A)=\left\{
\varphi\in C_c^1(A;\R^{m}):
\displaystyle  \max_{(x,y)\in A} |\varphi(x,y)|  \leq 1\right\}.
\]

\subsection{Relation between  $H$-perimeter and $\alpha$-perimeter}

Let $X_1,\ldots,X_h$ be the generators of an $H$-type Lie algebra, thought of
as left-invariant vector fields in $\Rn$ as in \eqref{HTYP}.
For an open  set $E\subset \R^n$ with Lipschitz boundary,
the Euclidean outer unit normal $N^E:\partial E\to\R^n$
is defined at $\mathcal H^{n-1}$-a.e.~point of $\partial E$.
We define the mapping $N_H^E: \partial E\to\R^h$
\[
       N_H ^E= (\langle N^E, X_1\rangle,\ldots,\langle N^E,X_h\rangle  ).
\]
Here, $\langle\cdot,\cdot\rangle$ is the standard scalar product
of $\Rn$ and $X_i$ is thought of as an element of $\Rn$ with respect
to the standard basis $\partial_{1},\dots,\partial_n$.

\begin{prop}\label{HH}
If  $E\subset\R^{n}$ is a bounded  open set with Lipschitz boundary
then the $H$-perimeter of $E$ in $\R^n$ is
\begin{equation}
 \label{HHH}
  \PH(E)=\int_{\partial E} |N_H^E (x,y)| \, d\mathcal H^{n-1},
\end{equation}
where   $\mathcal H ^{n-1}$  is the standard
$(n-1)$-dimensional Hausdorff measure in $\R^{n}$.
\end{prop}

\proof

The proof of \eqref{HHH} is standard and we only sketch
it.
The inequality
\[
 \PH(E) \leq \int_{\partial E} |N_H ^E(x,y)| \, d\mathcal H^{n-1}
\]
follows by the Cauchy-Schwarz inequality applied to the right hand side
of the identity
\[
\int_E \sum_{i=1}^h X_i \varphi_i\; dxdy = \int_{\partial E} \langle
N_H^E,\varphi \rangle d\mathcal H^{n-1},
\]
that holds for any $\varphi \in \F_h(\R^n)$.

The opposite inequality follows by
approximating $N_H^E/|N_H^E|$ with functions in $\F_h(\R^n)$.
In fact, by a  Lusin-type and Titze-extension
argument, for any $\varepsilon>0$
there exists $\varphi\in\F_h(\R^n)$ such that
\[
\int_{\partial E} \langle N_H^E,\varphi \rangle d\mathcal H^{n-1}\geq
 \int_{\partial E} |N_H^E (x,y)| \, d\mathcal H^{n-1}-\varepsilon .
\]
\endproof

The outer normal $N^E$ can be split  in the following way
\[
          N ^E= (N_x^E,N_y^E)\quad \textrm{with $N_x^E\in\R^h$ and
$N_y^E\in\R^k$.}
\]
For any $\alpha>0$, we call the mapping $N_\alpha^E: \partial E\to\R^n$
\begin{equation}\label{ALPHA}
       N_\alpha ^E= (N_x^E, |x|^\alpha N_y^E)
\end{equation}
the $\alpha$-normal to $\partial E$.
The same argument used to prove \eqref{HHH}
also shows that 
\begin{equation}
 \label{eq:per.rep}
  \P(E)=\int_{\partial E} |N_\alpha^E (x,y)| \, d\mathcal H^{n-1},
\end{equation}
for any set  $E\subset\R^{n}$   with Lipschitz boundary.

\begin{rem}
Formulas \eqref{HHH} and \eqref{eq:per.rep} hold also when $\partial E$ is $\mathcal{H}^{n-1}$-rectifiable.
\end{rem}

 \begin{prop}\label{PAPPA}
For any 
$x$-spherically symmetric set   $E\in\S $ there holds
$\PH(E) = \P(E)$ with $\alpha=1$. 
\end{prop}

\proof
By a standard approximation, using the results of \cite{FSSC}, it is sufficient to prove the claim for
smooth sets, e.g., for a bounded set $E\subset\Rn$ with Lipschitz boundary. By \eqref{HHH} and \eqref{eq:per.rep}, the claim $\PH(E)=\P(E)$ with $\alpha=1$ reads
\begin{equation}\label{lop}
 \PH(E) = \int_{\partial E} \sqrt{|N_x^E|^2+|x|^2 |N_y^E|^2} d\mathcal H^{n-1},
\end{equation}
 where $N^E = (N_x^E,N_y^E) \in \R^h\times\R^k$ is the unit Euclidean normal to
$\partial E$.
By the representation formula \eqref{HHH}, we have  
\[
     \PH(E)  = \int_{\partial E} \Big( \sum_{i=1}^h \langle X_i, N^E\rangle
^2\Big)^{1/2} d\mathcal H^{n-1},
\]
where, by \eqref{HTYP},   for any $i=1,\dots,h$
\[ 
\begin{split}
\langle X_i, N^E \rangle^2 &  
 =  
\Big(  N_{x_i} ^E - \sum _{\ell=1}^k
\sum_{j=1}^h Q^\ell_{ij} x_j N_{y_\ell}^E\Big)^2
\\
&
=(N_{x_i} ^E)^2 - 2 N_{x_i} ^E
\sum _{\ell=1}^k
\sum_{j=1}^h Q^\ell_{ij} x_j N_{y_\ell}^E
+\Big(
\sum _{\ell=1}^k
\sum_{j=1}^h Q^\ell_{ij} x_j N_{y_\ell}^E\Big)^2,
\end{split}
\]
and thus
\begin{equation}\label{pappa}
\begin{split}
\sum_{i=1}^h \langle X_i, N^E \rangle^2 &  
 =  
| N_{x} ^E|^2   - 2 
\sum _{\ell=1}^k
\sum_{i,j=1}^h Q^\ell_{ij} x_j N_{x_i} ^EN_{y_\ell}^E
+ \sum_{i=1}^h
\sum _{\ell,m=1}^k
\sum_{j,p=1}^h 
  Q^\ell_{ij}Q^{m}_{ i  p} x_j 
   x_p N_{y_\ell}^E N_{y_m}^E.
\end{split}
\end{equation}

Since the set $E$ is $x$-spherically symmetric, the component
$N_x^E$ of the normal satisfies
the identity
\begin{equation}\label{NUMME}
    N_x^E =\frac{x}{|x|} |N_x^E|.
\end{equation}
The bilinear form $Q:\R^h\times\R^h\to\R^k$  is skew-symmetric, i.e., 
we have $Q(x,x') = -Q(x',x)$ for all $x,x'\in\R^h$ or, equivalently,
$Q_{ij}^\ell=-Q_{ji}^\ell$.
Using \eqref{NUMME}, it follows that for any $\ell=1,\ldots,k$ we have
\begin{equation}\label{pix1}
\sum_{i,j=1}^h Q^\ell_{ij} x_j N_{x_i} ^E =
\frac{|N_x^E|}{|x|} \sum_{i,j=1}^h Q^\ell_{ij} x_i x_j=0.
\end{equation}

Next, we insert into identity
\eqref{H-type}, that defines an $H$-type group, the vector fields 
\[
 X = X'= \sum_{i=1}^h x_i X_i,\quad 
 Y = \sum_{j=1}^k N_{y_j}^E Y_j,
\]
where $x \in\R^h$, $N_y^E =(N_{y_1}^E,\ldots, N_{y_k}^E)$, and
$X_i$, $Y_j$ are the orthonormal vector fields in \eqref{HTYP}. 
After some computations that are omitted, using the definition \eqref{Kaplan} of the Kaplan mapping,
we obtain the identity  
\begin{equation}
   \label{tipoH2}
   \sum_{\ell ,m=1}^k\sum_{i,j,p=1}^h 
Q^\ell_{ ij} Q^m_{i p } N_{y_\ell}^E N_{y_m}^E x_j x_p  =
|x|^2 |N_y^E|^2.
\end{equation}

From \eqref{pappa}, \eqref{pix1}, and \eqref{tipoH2} we deduce that
\[
\sum_{i=1}^h \langle X_i, N^E \rangle^2 =  
| N_{x} ^E|^2  +|x|^2 |N_y^E|^2 ,
\]
 and formula \eqref{lop} follows.
\endproof

\subsection{$\alpha$-Perimeter for symmetric sets}
Thanks to Proposition \ref{PAPPA}, from now on we
will consider only $\alpha$-perimeter.

We say that a set  $E\subset\Rn =\R^h\times\R^k$ 
is $x$- and $y$-spherically symmetric if there exists
a set $G\subset \R^+\times \R^+$ such that
\[
         E = \big\{(x,y)\in \R^n : (|x|,|y|) \in G\big\}.
\]
We call $G$ the generating set of $E$. In the following we will use the constant
\[
c_{hk}=hk\omega_h\omega_k,
\]
where $\omega_m = \L^m(\{x\in\R^m:|x|<1\})$, for   $m\in\N$.

\begin{prop}\label{thm:N->k+1}
Let $E\subset \R^n$ be a bounded  open set with
finite $\alpha$-perimeter
that is  $x$- and $y$-spherically symmetric
with generating set $G\subset \R^+\times\R^+$.
Then we have:
\begin{equation}
\label{eq:N->2preliminary}
 \P(E) = c_{hk}
 \sup_{\psi \in \F_{2}(\R^+\times\R^+)}
\int _G  \Big( s^{k-1}  \partial _ r
 \big(   r^{h-1}\psi _1\big)
+ r^{h-1+\alpha} \partial_s\big( s^{k-1}   \psi _{2}\big) \Big) \,
 drds.
\end{equation}
In particular, if $E$ has Lipschitz boundary then we
have:
\begin{equation}
\label{eq:N->2}
\P(E) = c_{hk}
   \int_{\partial G} |(N_r^G,r^\alpha
N_s^G)|r^{h-1}s^{k-1}\;d\mathcal H^1(r,s),
\end{equation}
where  $N^G=(N_r^{G},N_s^{G})\in \R^{2}$ is the outer unit normal
to the boundary $\partial G\subset \R^+\times \R^+$.
\end{prop}

\proof
We prove a preliminary version of \eqref{eq:N->2preliminary}.
We claim that if $E$ is of finite $\alpha$-perimeter
and  $x$-spherically symmetric with generating set $F\subset \R^+\times\R^k$, then
we have:
\begin{equation}
   \label{eq:N->k+1}
  \P (E) =  h\omega_h
 \sup_{\psi \in \F_{1+k}(\R^+\times\R^k)}
\int _F  \Big( \partial_r\big( r^{h-1} \psi _1\big)
+ r^{h-1+\alpha} \sum_{j=1}^k \partial_{y_j}  \psi_{1+j} \Big)  drdy
  =\Q(F),
\end{equation}
where $\Q$ is defined via the last identity.
For any test function
$\psi \in\F_{1+k} (\R^+\times\R^k)$ we define
the test function $\varphi\in \F_n(\R^n)$
\begin{equation}
 \label{eq:psiphi}
 \varphi(x,y)=\left(\frac{x}{|x|}\psi_1(|x|,y),\psi_2(|x|,y),\dots,
 \psi_{1+k}(|x|,y)\right)\text{ for }|x|\neq0,
 \end{equation}
and $\varphi(0,y)=0$.
For any $i=1,\ldots,h$, $j=1,\ldots,k$, and $x\neq 0$, we have the identities
\begin{align*}
  &\partial_{x_i}\varphi_{i}(x,y)=
\left(\frac{1}{|x|}-\frac{x_i^2}{|x|^3}\right)\psi_1(|x|,y)+\frac{x_i^2}
{|x|^2}\partial_r\psi_1(|x|,y),
\\
&\partial_{y_j}\varphi_{h+j}(x,y)
=\partial_{y_j}\psi_{1+j}(|x|,y),
\end{align*}
and thus, the $\alpha$-divergence defined by
\begin{equation}\label{DIVA}
\diva\varphi(x,y)=\sum_{i=1}^h \frac{\partial \varphi_i (x,y)}{\partial x_i} +|x|^\alpha \sum_{j=1}^k \frac{\partial \varphi_{h+j}(x,y)}{\partial y_j}
\end{equation}
satisfies
\begin{equation}
\label{eq:div}
\diva\varphi(x,y)=
\dfrac{h-1}{|x|}\psi_1(|x|,y)+\partial_r\psi_1(|x|,
y)+|x|^\alpha\sum_{j=1}^k\partial_{y_j}\psi_{1+j}(|x|,y).
\end{equation}

For any $y\in\R^k$ we define the section $F^y = \big\{ r>0: (r,y)\in F\big\}$.
Using Fubini-Tonelli theorem, spherical coordinates
in $\R^h$, the symmetry of $E$,
and  \eqref{eq:div} we obtain
\begin{equation}
\label{eq:intdiv}
\begin{split}
\int_E \diva\varphi\;dxdy
&=\int_{\R^k}\int_{F^y} \int_{|x|=r}
\left(
\dfrac{h-1}{r}
\psi_1+\partial_r\psi_1+r^\alpha\sum_{j=1}^k
\partial_{y_j
}\psi_{1+j}\right)\;d\HH^{h-1}(x)dr dy
\\
&=h\omega_h\int_{\R^k}\int_{F^y}r^{h-1}
\left(
\dfrac{h-1}{r}\psi_1+\partial_r\psi_1+r^\alpha
\sum_{j=1}^k\partial_{y_j}\psi_{1+j}\right)\;dr dy
\\
&=h\omega_h\int_{F}\partial_r(r^{h-1}\psi_1)+r^{\alpha+h-1}\sum_{j=1}
^k\partial_{y_j}\psi_{1+j}\;dr dy.
\end{split}
\end{equation}
Because $\psi$ is arbitrary, this proves the inequality $\geq$
in \eqref{eq:N->k+1}.

We prove the opposite inequality when $E\subset\Rn$
is an $x$-symmetric bounded open set with smooth boundary.
The unit outer normal $N^E=(N_x^E,N_y^E)$ is continuously
defined on
$\partial E$. At points $(0,y)\in \partial E$, however,
 we have $N_x^E(0,y)=0$ and thus $N_\alpha^E(0,y)=0$.
For any  $\varepsilon>0$ we consider the compact set
 $K= \big\{ (x,y) \in \partial E : |x| \geq \delta\big\}$,
where $\delta>0$ is such that $\P(E; \{|x|=\delta \})=0$
and
\begin{equation} \label{EPS}
\int_{\partial E\setminus K}
|N_\alpha^E(x,y)| \, d\mathcal H^{n-1}<\varepsilon.
\end{equation}

Let $H\subset \R^+\times\R^k$ be the generating set of $K$.
By standard extension theorems, there exists
$\psi \in \F_{1+k}(\R^+\times\R^k)$
such that
\[
 \psi(r,y) = \frac{(N_r^F(r,y), r^\alpha N_y^F(r,y))}
 {|(N_r^F(r,y), r^\alpha N_y^F(r,y)|} \quad \textrm{for } (r,y)\in H.
\]
The mapping   $\varphi \in\F_n(\Rn)$
introduced in \eqref{eq:psiphi} satisfies
\begin{equation}\label{fia}
 \varphi(x,y) = \frac{N_\alpha^E(x,y)}{|N_\alpha^E(x,y)|},\quad \textrm{for }(x,y)\in K.
\end{equation}

Then, by identity \eqref{eq:intdiv}, the divergence theorem,
\eqref{fia},  \eqref{EPS}, and \eqref{eq:per.rep}
 we have
\[
\begin{split}
 \Q(F)  & \geq
\int _F  \Big( \partial_r\big( r^{h-1} \psi _1\big)
+ r^{h-1+\alpha} \sum_{j=1}^k \partial_{y_j}  \psi_{1+j} \Big)  drdy
\\
& = \int_E \diva \varphi \, dx dy  = \int _{\partial E} \langle\varphi , N_\alpha^E\rangle \, d\mathcal
H^{n-1}
\\
&
= \int_{K}|N_\alpha^E(x,y)|\,
  d\mathcal H^{n-1}
 +\int _{\partial E \setminus K }\langle \varphi, N_\alpha^E\rangle \,
 d\mathcal H^{n-1}
\\
&
  \geq \P(E)- 2\varepsilon.
\end{split}
\]
This proves \eqref{eq:N->k+1} when $\partial E$ is smooth.
The general case follows by approximation.
Let $E\subset\Rn$ be a set of finite
$\alpha$-perimeter and finite Lebesgue measure that is $x$-symmetric
with generating set $F\subset\R^+\times\R^k$.
By \cite[Theorem 2.2.2]{FSSC}, there exists
a sequence $(E_j)_{j\in\N}$ such that each $E_j$ is of class $C^\infty$
\[
 \lim_{j\to\infty} \L^n(E_j\Delta E) = 0 \quad \textrm{and}\quad
  \lim_{j\to\infty} \P(E_j) = \P(E).
\]
Each $E_j$ can be also assumed to be $x$-spherically symmetric with generating
set $F_j\subset\R^+\times\R^k$. Then we also have
\[
       \lim_{j\to\infty} \mathcal L^{1+k}(F_j\Delta F) = 0 .
\]
By lower semicontinuity and \eqref{eq:N->k+1} for the smooth case, we have
\[
              \Q(F) \leq \liminf_{j\to\infty}
                          \Q(F_j)
                       =\lim_{j\to\infty} \P(E_j)=\P(E).
\]
This concludes the proof of \eqref{eq:N->k+1} for any set
$E$ with finite $\alpha$-perimeter.

The general formula \eqref{eq:N->2preliminary} for sets that
are also $y$-spherically symmetric
can be proved in a similar way and we can omit the details.

Formula \eqref{eq:N->2} for sets $E$ with Lipschitz boundary
follows from \eqref{eq:N->2preliminary} with the same argument sketched
in the proof of Proposition \ref{eq:per.rep}.
\endproof

\subsection{$\alpha$-Perimeter in the case $h=1$}

When $h=1$ there exists a change of coordinates that transforms
$\alpha$-perimeter into the standard perimeter (see
\cite{MontiMorbidelli} for the case of the plane $h=k=1$). 
Let $n=1+k$ and consider
the mappings $\Phi,\Psi:\Rn\to\Rn$
\begin{eqnarray*}
\Psi(x,y)=\left(\mathrm{sgn}(x)\dfrac{|x|^{\alpha+1}}{\alpha+1},y\right)
\quad \text{and}\quad
\Phi(\xi,\eta)=\left(\mathrm{sgn}(\xi)|(\alpha+1)\xi|^{\frac{1}{\alpha+1}},
\eta\right).
\end{eqnarray*}
Then we have   $\Phi\circ\Psi=\Psi\circ\Phi=\mathrm{Id}_{\R^{n}}$.

\begin{prop}\label{prop:alpha-euclidean}
Let $h=1$ and $n=1+k$.
For any measurable set   $E\subset\R^{n}$ we have
\begin{equation}\label{P=P}
 \P (E)=  \sup \Big\{ \int_{\Psi(E)} \mathrm{div} \psi\, d\xi d\eta\, :\,
 \psi\in\F_{n}(\Rn)\Big\}.
\end{equation}
\end{prop}
\proof
First notice that the supremum in the right hand side can be equivalently
computed over all vector fields $\psi:\Rn\to\Rn$
in the Sobolev space $W^{1,1}_0(\Rn;\Rn)$ such that $\| \psi\|_\infty\leq 1$.

 For any $\varphi \in \F_n(\Rn)$, let $\psi = \varphi\circ\Phi$.
Then for any $j=1,\ldots, k = n-1$, we have
\begin{equation}
\begin{split}
 \partial_{\xi} \psi_1 (\xi, \eta)&
=\partial_{\xi}\big( \varphi_1\circ\Phi\big) (\xi,\eta)
  =
|(\alpha+1)\xi|^{-\frac{\alpha}{\alpha+1}}\partial_x\varphi_1(\Phi(\xi,\eta)),
\\
\partial_{\eta_j}\psi_{1+j} (\xi, \eta)& =
\partial_{\eta_j} \big(\varphi_{1+j}\circ \Phi\big)(\xi,\eta)
 =\partial_{y_j}\varphi_{1+j}(\Phi(\xi,\eta)).
\end{split}
\end{equation}
In particular, we have $\psi \in W^{1,1}_0(\Rn;\Rn)$ and
$\|\psi\|_\infty\leq 1$. Then, the standard divergence of $\psi$
satisfies
\[
 \mathrm{div}\psi(\xi,\eta) = |(\alpha+1)\xi|^{-\frac{\alpha}{\alpha+1}} \diva \phi (\Phi(\xi,\eta)).
\]
The determinant Jacobian of the change of variable
$(x,y)=\Phi(\xi,\eta)$ is
\begin{equation}\label{detJ}
|\det J\Phi(\xi,\eta)| =
|(\alpha+1)\xi| ^{-\frac{\alpha}{\alpha+1}}.
\end{equation}
and thus
we obtain
\begin{equation}
 \begin{split}
 \int_{E} \diva\varphi(x,y)\;dxdy
  & =\int_{\Psi(E)} \diva\varphi(\Phi(\xi,\eta))
|\det J\Phi(\xi,\eta)|\, d\xi d\eta
\\
& =
 \int_{\Psi(E)} \mathrm{div}\psi (\xi,\eta) \;d\xi d\eta.
\end{split}
\end{equation}
The claim follows.

\endproof

\section{Rearrangements}\label{TRE}

In this section, we prove various rearrangement inequalities for $\alpha$-perimeter in $\Rn$. 
We consider first the case $h=1$.
In this case, there are
a Steiner type rearrangement in the $x$-variable
and a Schwartz rearrangement in the $y$ variables
that reduce the isoperimetric problem in $\Rn$ to
a problem  for Lipschitz graphs in the first quadrant $\R^+\times\R^+$.
Then we consider dimensions $h\geq 2$, where we can  rearrange  
sets in $\R^h$   that are  already
$x$-spherically symmetric.

\subsection{Rearrangement in the case $h=1$}
Let $h=1$ and $n=1+k$.
We say that a set
$E\subset\Rn$   is $x$-symmetric if $(x,y)\in E$ implies $(-x,y)\in E$;
we say that $E$ is  $x$-convex if the section
 $E^{y}=\{x\in\R : (x,y)\in E\}$ is an interval for every
$y\in\R^k$; finally, we say that $E$ is $y$-Schwartz symmetric
if  for every $x\in\R$ the section $E^x=\{y\in\R^{k} :
(x,y)\in E\}$ is an (open)
Euclidean ball in $\R^k$ centered at the origin.

\begin{thm}\label{thm:rearrangement}
Let $h=1$ and $n=1+k$. For any set  $E\subset\R^{n}$
such that $\P(E)<\infty$ and
$0<\L^n(E)<\infty$ there exists an $x$-symmetric, $x$-convex, and
$y$-Schwartz symmetric set   $E^*\subset\R^{n}$ such that $\P(E^*)\leq
\P(E)$ and $\L^n(E^*)=\L^n(E)$.

Moreover, if $\P(E^*)=
\P(E)$ then
$E$ is $x$-symmetric, $x$-convex and
there exist functions $c:[0,\infty)\to\R^k$ and $f:[0,\infty)\to[0,\infty]$
such that for $\L^1$-a.e.~$x\in \R$ we have 
\begin{equation}
 \label{sezsez}
 E^ x = \{ y \in \R^k: |y-c(|x|)|< f(|x|)\}.
\end{equation}
\end{thm}

\proof

By Proposition \ref{prop:alpha-euclidean},
the set $F=\Psi(E)\subset\Rn$ satisfies
$P(F)=\P(E)$, where $P$ stands for the standard perimeter in $\Rn$.
We define the measure $\mu$ on $\R^{n}$
\begin{equation}\label{mu-vol}
       \mu(F)=\int_F|(\alpha+1)\xi|^{-\frac{\alpha}{\alpha+1}}\,d\xi d\eta.
\end{equation}
Then, by \eqref{detJ} we also have the identity
 $\mu(F)=\L^n(E)$.

We rearrange the set $F$ using Steiner symmetrization in direction $\xi$.
Namely, we let
\[
  F_1= \{(\xi,\eta)\in\R^n : |\xi|<\mathcal{L}^1(F^\eta)/2\},
\]
where $F^\eta = \{\xi \in\R: (\xi,\eta)\in F\}$.
The set $F_1$ is $\xi$-symmetric and $\xi$-convex.
By classical results on Steiner symmetrization
we have  $P(F_1)\leq P(F)$ and the equality $P(F_1)= P(F)$
implies that $F$ is $\xi$-convex: namely, a.e.~section $F^\eta$ is (equivalent
to) an interval.

The $\mu$-volume of $F_1$ is
\[
     \mu(F_1)=\int_{F_1}|(\alpha+1)\xi|^{-\frac{\alpha}{\alpha+1}}
                  d\xi d\eta
       =\int_{\R^k}\Big(\int_{F_1^\eta}
 |(\alpha+1)\xi|^{-\frac{ \alpha}{\alpha+1}}d\xi\Big)d\eta.
\]
For any  measurable set $I\subset\R$ with
finite measure, the symmetrized set $I^*=(-\L^1(I)/2,\L^1(I)/2)$
satisfies the following inequality
(see
\cite{MontiMorbidelli}, page 361)
\begin{equation}
\label{mu-int}
 \int_{I}|\xi|^{-\frac{\alpha}{\alpha+1}} d\xi\leq
 \int_{I^*}|\xi|^{-\frac{\alpha}{ \alpha+1}} d\xi.
\end{equation}
Moreover, if $\L^1(I\Delta I^*)>0$ then the inequality is strict.
This implies that   $\mu(F_1)\geq\mu(F)$
and the inequality is strict if $F$ is not equivalent to an $\xi$-symmetric
and $\xi$-convex set.

We rearrange the set $F_1$ using Schwartz symmetrization
in $\R^k$, namely we let
\[
 F_2=\Big\{
 (\xi,\eta)\in\R^{n} :
|\eta|<\Big(\frac{\L^k(F^\xi_1)}{\omega_k}\Big)^{\frac{1}{k}}\Big\}.
\]
By classical results on Schwartz rearrangement, we have
$P(F_2)\leq P(F_1)$ and the equality
$P(F_2)= P(F_1)$ implies that a.e.~section $F_1^\xi$ is an Euclidean ball
\begin{equation}\label{section1}
      F_1^\xi= \{\eta \in\R^k: |\eta -d(|\xi|)|<\varrho(|\xi|)\}
\end{equation}
for some $d(|\xi|)\in\R^k$ and $\varrho(|\xi|)\in[0,\infty]$.
By Fubini-Tonelli theorem, the $\mu$-volume is preserved:
\begin{equation}
       \mu(F_2)
 =\int_{\R}|(\alpha+1)\xi|^{-\frac{\alpha}{\alpha+1}}
\L^k(F_2^{\xi}) d\xi
=\int_{\R}|(\alpha+1)\xi|^{-\frac{\alpha}{\alpha+1}}
\L^k(F_1^{
\xi})d\xi=\mu(F_1).
\end{equation}

Recall that $\delta_\lambda(x,y)=(\lambda x,\lambda^{\alpha+1}y)$.
The set
$E^*=\delta_\lambda(\Phi(F_2))$, with $\lambda>0$ such that
  $\L^n(E^*)=\L^n(E)$, satisfies the claims in the statement of the theorem.
In fact, we have   $0<\lambda\leq1$ because
\[
\L^n( \Phi(F_2))
=\mu(F_2)=\mu(F_1)\geq \mu(F)=\L^n(E),
\]
and then, by the scaling property of $\alpha$-perimeter  we have
\[
\P (E^*)
=\lambda^{d-1}
\P (\Phi(F_2)) \leq \P(\Phi(F_2)) = P(F_2)\leq P(F_1)\leq P(F)=\P (E).
\]
This proves the first part of the theorem.

If $\P(E^*) = \P(E)$ then we have $P(F_2) = P(F_1)$ and $\lambda=1$.
From the first equality we deduce that the sections $F_1^\xi$ are
of the form \eqref{section1} and claim \eqref{sezsez} holds with
$c(|x|)=d\big(|x|^{\alpha+1} / (\alpha+1) \big)$
and $f(|x|) = \varrho\big(|x|^{\alpha+1} / (\alpha+1) \big)$.
From $\lambda=1$ we deduce that
\[
\mu(F) = \L^n(E) = \L^n(E^*) = \L^n(\Phi(F_2)) = \mu(F_2) = \mu(F_1),
\]
and thus $F$ is $\xi$-symmetric and $\xi$-convex.
The same holds then for $E$.

\endproof

\subsection{Rearrangement in the case $h\geq2$}
We prove the analogous of Theorem \ref{thm:rearrangement}
when $h\geq 2$. We need to start from a set $E\subset\Rn$ that is
$x$-spherically symmetric
\[
        E = \{(x,y)\in\Rn : (|x|, y) \in F\}
\]
for some generating set $F\subset \R^+\times \R^k$.

By the proof of Proposition \ref{thm:N->k+1}, see \eqref{eq:N->k+1},
we have the identity $\P(E) = \Q(F)$, where
\begin{equation} \label{ququ}
\Q(F) =
h\omega_h
 \sup_{\psi \in \F_{1+k}(\R^+\times\R^k)}
\int _F  \Big( \partial_r\big( r^{h-1} \psi _1\big)
+ r^{h-1+\alpha} \sum_{j=1}^k \partial_{y_j}  \psi_{1+j} \Big)  drdy.
\end{equation}

Our goal is to improve the $x$-spherical symmetry
to the $x$-Schwartz symmetry.
A set $E\subset\Rn$ is $x$-Schwartz symmetric
if for all $y\in\R^k$ we have
\[
        E^y =\{x\in \R^h: (x,y)\in E\} = \{ x\in \R^h: |x|<\varrho(y)\}
\]
for some function $\varrho:\R^k\to[0,\infty]$.
To obtain the Schwartz symmetry, we use the radial rearrangement
technique introduced in \cite{M3}.

\begin{thm}\label{thm:rear2}
Let $h\geq 2$, $k\geq1$ and $n=h+k$. For any set  $E\subset\R^{n}$
that is $x$-spherically symmetric and such that $\P(E)<\infty$ and
$0<\L^n(E)<\infty$ there exists an $x$-  and $y$-Schwartz symmetric
set $E^*\subset\R^{n}$ such that $\P(E^*)\leq
\P(E)$ and $\L^n(E^*)=\L^n(E)$.

Moreover, if $\P(E^*)=
\P(E)$ then
$E$ is $x$-Schwartz symmetric  and
there exist functions $c:[0,\infty)\to\R^k$ and $f:[0,\infty)\to[0,\infty]$
such that, up to a negligible set,
we have
\begin{equation}\label{TOPOLINO}
  E = \{ (x,y) \in \Rn: |y-c(|x|)|< f(|x|)\}.
\end{equation}
\end{thm}

\proof
Let $F\subset \R^+\times\R$ be the generating set of $E$.
We define the volume of $F$ via the following formula
\[
 V(F) = \omega_h \int _F r^{h-1} dr dy = \mathcal L^{n}(E).
\]

We rearrange $F$ in the coordinate $r$ using the linear density
$r^{h-1+\alpha}$ that appears, in \eqref{ququ},
in the part of divergence depending on the coordinates $y$.
Namely, we define the function $g:\R^k\to[0,\infty]$ via the identity
\begin{equation} \label{fifa}
    \frac{1}{h+\alpha} g(y)^{h+\alpha} =   \int _0^{g(y)} r^{h-1+\alpha} dr   =
\int _{F_y}r^{h-1+\alpha} dr,
\end{equation}
and we let
\[
            F^\sharp =\big\{ (r,y)\in \R^+\times\R^k : 0<r< g(y)\big\}.
\]

We claim that $\Q(F^\sharp)\leq \Q(F)$ and $V(F^\sharp)\geq V(F)$, with equality
$V(F^\sharp) = V(F)$ holding if and only if $F^\sharp=F$, up to a negligible set.

For any open set $A\subset \R^+\times\R^k$, we define
\begin{equation}
\label{perpar}
\begin{split}
  \Q_0(F;A) &
 =  \sup_{\psi \in \F_{1}(A)}
\int _F    \partial_r\big( r^{h-1} \psi\big)
\,   drdy,
\\
\Q_{j}(F;A) &
  =\sup_{\psi \in \F_{1}(A)}
\int _F
 r^{h-1+\alpha} \partial_{y_j}  \psi\, drdy, \quad j=1,\ldots,k.
\end{split}
\end{equation}
The open sets mappings $A\mapsto
\Q_{j}(F;A)$, $j=0,1,\ldots,k$, extend to Borel measures.
For any Borel set $B\subset\R^k$
and $j=0,1,\ldots,k$, we define the measures
\[
\begin{split}
 \mu_j(B) & = \Q_j(F;\R^+\times B ),
\\
 \mu_j^\sharp(B) & = \Q_j(F^\sharp;\R^+\times B ).
\end{split}
\]
By Step 1 and Step 2 of the proof of Theorem 1.5  in \cite{M3}, see page 106, we have
$\mu_j^\sharp(B) \leq \mu_j(B)$ for any Borel set $B\subset \R^k$
and for any $j=0,1,\ldots,k$.
It follows that
 the vector valued Borel measures $\mu = (\mu_0,\ldots,\mu_k)$
and $\mu^\sharp = (\mu_0^\sharp,\ldots,\mu_k^\sharp)$ satisfy
\[
   |\mu^\sharp|(\R^k)\leq |\mu|(\R^k),
\]
where $|\cdot|$ denotes the total variation.
This is equivalent to $\Q(F^\sharp)\leq \Q(F)$.

We claim that for any $y\in\R^k$ we have
\begin{equation}
 \label{sphinx}
\frac 1 h g(y)^h=   \int_{F^\sharp_y} r^{h-1} \, dr \geq \int_{F_y} r^{h-1} \,
dr,
\end{equation}
with strict inequality unless $F^\sharp_y = F_y$ up to a negligible set.
From \eqref{sphinx}, by Fubini-Tonelli theorem it follows
that
$V(F^\sharp)\geq V(F)$ with strict inequality
unless $F^\sharp=F$ up to a negligible set.
By \eqref{fifa}, claim \eqref{sphinx} is equivalent to
\begin{equation}
\label{top}
   \Big((h+\alpha)      \int_{F_y} r^{h-1+\alpha} dr\Big) ^{\frac{1}{h+\alpha}}
 \geq \Big(h     \int_{F_y} r^{h-1} dr\Big) ^{\frac{1}{h}},
\end{equation}
and this inequality holds for any measurable set $F_y\subset\R^+$,
for any $h\geq2$,  and $\alpha>0$, by
Example 2.5 in \cite{M3}. Moreover, we have equality
in \eqref{top} if and only if $F_y=(0,g(y))$.

Let $E_1^\sharp\subset\R^n$  be
the $x$-Schwartz symmetric set  with generating set
$F^\sharp$. Then we have
\[
\L^n(E _1^\sharp) = V(F^\sharp)\geq V(F) = \L^n(E),
\]
with strict inequality unless $F^\sharp= F$.
Then there exists $0<\lambda\leq 1$ such that the set
$E^\sharp =\delta_\lambda (E_1^\sharp)$ satisfies
$\L^n(E^\sharp) = \L^n(E)$. Since $\lambda \leq 1$, we also have
\[
  \P(E^\sharp) = \lambda^{d-1} \P(E_1^\sharp) \leq \P(E_1^\sharp)
                 = \Q(F^\sharp) \leq \Q(F) = \P(E).
\]
If $\P(E^\sharp) = \P(E)$ then it must be $\lambda=1$ and thus
$F^\sharp = F$, that in turn implies $E^\sharp= E$, up to a negligible set.

Now the theorem can be concluded applying to
$E^\sharp$ a Schwartz
rearrangement in the variable $y\in\R^k$. This rearrangement is standard,
see the general argument in \cite{M4}. 
The
resulting set
$E^*\subset\Rn$ satisfies $\P(E^*)\leq \P(E)$ and
also the other claims in the theorem. 

\endproof

\section{Existence of isoperimetric sets}
\label{sec:existence}

In this section, we prove   existence
of solutions to the isoperimetric problem for $\alpha$-perimeter and
$H$-perimeter.
When $h\geq 2$, we prove the existence of solutions in the class
of $x$-spherically symmetric sets.
The proof is based on a concentration-compactness
argument.

For any set $E\subset \Rn $ and $t>0$, we let
\begin{equation}\begin{split}
 & E_{t-}^{x}=\left\{(x,y) \in E : |x|<t\right\} 
  \quad\text{and}\quad 
  E_{t}^{x}=\left\{(x,y) \in E :
|x|=t\right\},
\\& 
E_{t-}^{y}=\left\{(x,y) \in E : |y|<t\right\}
 \quad\text{and}\quad 
E_{t}^{y}=\left\{(x,y) \in E : |y|=t\right\}.
\end{split}
\end{equation}
We also define  
\begin{equation} \label{paperinik2}
  v^x_E(t) = \mathcal{H}^{n-1}(E_t^x ),
\end{equation}
and
\begin{equation}\label{paperinik3}
  v^y_E(t) =  \int _{E_t^y } |x|^\alpha d \mathcal
H^{n-1}.
\end{equation}
In the following, we use the short notation $\{|x|<t\} = \{(x,y)\in\Rn:|x|<t\}$
and $\{|y|<t\} = \{(x,y)\in\Rn:|y|<t\}$.

\begin{prop}
\label{prop:derivatives}
Let $E\subset \Rn$ be a set with finite measure and finite $\alpha$-perimeter.
Then for a.e.~$t>0$
we have
\begin{equation}
 \label{eq:density}
 \begin{split}
 \P(E_{t-}^x )=
\P(E;E_{t-}^x )+ {v}^x_E(t)\quad \text{ and }\quad 
 \P(E_{t-}^y )=
\P(E;E_{t-}^y)+ {v}^y_E(t).
\end{split}
\end{equation}
\end{prop}

\proof
We prove the claim for $E^y_{t-}$.
Let $\{ \phi_\varepsilon\}_{\varepsilon>0}$ be a standard family of
mollifiers in $\Rn$ and let
\[
     f_\varepsilon(z) = \int_E \phi_\varepsilon (|z-w|) dw,\quad z\in\Rn.
\]
Then $f_\varepsilon\in C^\infty(\Rn)$ and
$f_\varepsilon\to \chi_E$ in $L^1(\Rn)$ for $\varepsilon\to0$.
Therefore, by the coarea formula we also have, for a.e.~$t>0$
and possibly for a suitable infinitesimal sequence of $\varepsilon$'s,
\begin{equation}
 \label{pippo1}
  \lim_{\varepsilon\to0} \int _{\{|y|=t\}} |f_\varepsilon
-\chi_{E}|
d\mathcal H^{n-1}=0.
\end{equation}

Since $E$ has finite $\alpha$-perimeter, the set
$\{t>0: \P(E; \{ |y|=t\})>0\}$ is at most countable, and thus
\begin{equation}
\label{pippo2}
\P(E;\{ |y|=t\})=0\quad \text{for a.e.~$t>0$.}
\end{equation}

We use the notation $\nabla\! _\alpha f_\varepsilon =
(X_1 f_\varepsilon ,\ldots,X_h f_\varepsilon, Y_1  f_\varepsilon,
\ldots Y_k f_\varepsilon)$, where $X_i$, $Y_j$ are the vector fields \eqref{GTYP}.
By the divergence Theorem, for any $\varphi\in C^1_c(\Rn,\Rn)$ we have
\begin{equation}
\label{topolino1}
\begin{split}
 \int_{\{ |y|<t\}  } f_\varepsilon(z) \diva \varphi(z)\, dz
  & =\int_{\{|y|<t\} }  \big( \diva( f_\varepsilon \varphi) -
\langle \nabla\!_\alpha f_\varepsilon, \varphi\rangle  \big) dz
 \\
&
= - \int _{\{ |y|=t\} }  f_\varepsilon(z) |x|^{\alpha} \langle N,\varphi (z)\rangle
   d\mathcal H^{n-1}
-
\int_{\{|y|<t\} }
\langle \nabla\!_\alpha f_\varepsilon, \varphi\rangle  dz,
\end{split}
\end{equation}
where $N=(0,-y/|y|)$ is the inner unit normal of $\{ |y|< t\}$.
For any $t>0$, we have
\begin{equation}
\label{topolino2}
 \lim_{\varepsilon\to0}
\int_{\{ |y| <t\}  } f_\varepsilon(z) \diva \varphi(z)\, dz
=
\int_{E^y_{t-} } \diva \varphi(z)\, dz ,
\end{equation}
and, for any $t>0$ satisfying \eqref{pippo1},
\begin{equation}
\label{topolino3}
 \lim_{\varepsilon\to0}
 \int _{\{ |y|=t\} }  f_\varepsilon(z) |x|^{\alpha} \langle N,\varphi (z)\rangle
   d\mathcal H^{n-1}
=
 \int _{E^y_t }  |x|^{\alpha} \langle N,\varphi (z)\rangle
   d\mathcal H^{n-1} .
\end{equation}

On the other hand, we claim that
\begin{equation}\label{pippo3}
\lim_{\varepsilon\to0}\int_{\{|y|<t\} }
\langle \nabla\!_\alpha f_\varepsilon, \varphi\rangle  dz=
\int_{\{|y|<t\} } \Big\{
\sum_{i=1}^ h \varphi _i d\mu_E^{x_i} + \sum_{\ell=1}^k \varphi_{h+\ell} |x|^\alpha
d\mu_E^{y_\ell}\Big\},
\end{equation}
where $\mu_E^{x_i}$ and $\mu_{E}^{y_\ell}$ are the distributional partial derivatives
of $\chi_E$, that are Borel measures on $\Rn$,
because $E$ has finite $\alpha$-perimeter.
For the coordinate  $y_\ell$, we have
\[
\begin{split}
 \int_{\{|y|<t\}} \varphi_{h+\ell}(z)
|x|^\alpha \partial _{y_\ell} f_\varepsilon(z)
dz
& =
 \int_{\{|y|<t\}} \varphi_{h+\ell}(z)
|x|^\alpha \int _{E} \partial_{y_\ell} \phi_\varepsilon(|z-w|)dw\, dz
\\
&
=-\int_{\{|y|<t\}} \varphi_{h+\ell}(z)
|x|^\alpha \int _{E} \partial_{\eta_\ell} \phi_\varepsilon(|z-w|)dw\, dz
\\&
=\int_{\{|y|<t\}} \varphi_{h+\ell}(z)
|x|^\alpha \int _{\Rn } \phi_\varepsilon(|z-w|)d\mu_E^{y_\ell}(w)
\, dz
\\&
=\int _{\Rn } \int_{\{|y|<t\}} \varphi_{h+\ell}(z)
|x|^\alpha
\phi_\varepsilon(|z-w|)dz\, d\mu_E^{y_\ell}(w),
\end{split}
\]
where we let $w=(\xi,\eta)\in\R^h\times\R^k$.
By \eqref{pippo2}, the measure $\mu_E^{y_\ell}$
is concentrated on $\{ |y|\neq t\}$. It follows that
\[
\lim_{\varepsilon\to 0}
\int _{\Rn } \int_{\{|y|<t\}} \varphi_{h+\ell}(z)
|x|^\alpha
\phi_\varepsilon(|z-w|)dz\, d\mu_E^{y_\ell}(w)
=
\int_{\{|\eta|<t\}} \varphi_{h+\ell}(w)
|\xi|^\alpha  d\mu_E^{y_\ell}(w).
\]
This proves \eqref{pippo3}.

Now, from \eqref{topolino1}--\eqref{pippo3} we deduce that
\begin{equation}\label{PPP}
\begin{split}
    \int_{E\cap \{ |y| <t\}  } \diva \varphi(z)\, dz &
 =
-\int _{E\cap \{ |y|=t\} }  |x|^{\alpha} \langle N,\varphi (z)\rangle
   d\mathcal H^{n-1}
\\
&
\qquad
-
\int_{\{|y|<t\} } \Big\{
\sum_{i=1}^ h \varphi _i d\mu_E^{x_i} + |x|^\alpha
\sum_{\ell=1}^k \varphi_{h+\ell}
d\mu_E^{y_\ell}\Big\},
\end{split}
\end{equation}
and the claim follows by optimizing the right hand side
over $\varphi\in \F_n(\Rn)$.

\endproof
\medskip

\begin{prop}\label{prop:calibration}
Let $E\subset\Rn$ be a  set
with finite measure and finite $\alpha$-perimeter.
For a.e.~$t>0$ we have   $\P(E^x_{t-})\leq \P(E)$ and  $\P(E^y_{t-})\leq
\P(E)$.
\end{prop}

\proof
The proof is a calibration argument. 
Notice that
\[
\begin{split}
\P (E^y_{t-})&=\P(E^y_{t-};\{|y|<t\})
+\P(E^y_{t-};\{|y|\geq t\})
\\
&=\P(E;\{|y|<t\})+\P(E^y_{t-};\{|y|=t\}).
\end{split}
\]
Let $t>0$ be such that $\P(E;\{|y|=t\})=0$; a.e.~$t>0$
has this property, see \eqref{pippo2}. It is sufficient to show that
\[
\P(E^y_{t-};\{|y|=t\})\leq\P(E;\{|y|\geq t\})=\P(E;\{|y|>t\}).
\]

The function   $\varphi(x,y) = (0,-y/|y|)\in\Rn$,  
$|y|\neq 0$, has negative divergence:
\[
 \diva\varphi(x,y)=-|x|^\alpha\sum_{\ell=1}^k\Big(\frac{1}{|y|}-\frac{y_{\ell
}^2}{|y|^3}\Big)=-\frac{(k-1)|x|^\alpha}{|y|}\leq 0.
\]
As in the  proof  of  \eqref{PPP}, we have
\[
\begin{split}
0\geq\int _{E\cap\{ |y|>t\} } \diva\varphi\, dz &= \int_{ E^y_t }
 |x|^\alpha d\mathcal H^{n-1}-\int_{\{|y|>t\}}|x|^\alpha
 \sum_{\ell=1}^k\varphi_{h+\ell} d\mu_E^{y_\ell}
\\
&
\geq
\int_{ E\cap\{ |y|=t\}}
 |x|^\alpha d\mathcal H^{n-1}-\P(E; \{|y|>t\}).
\end{split}
\]
By the  representation formula \eqref{eq:per.rep},
we obtain
\[
\P(E^y_{t-};\{|y|=t\})=\int_{ E^y_t}
 |x|^\alpha d\mathcal H^{n-1}\leq \P(E; \{|y|>t\}).
\]
This ends the proof.
\endproof

We prove the existence of isoperimetric sets
assuming   the validity of the following 
isoperimetric inequality, holding for any $\mathcal{L}^n$-measurable
set $E\subset\Rn$ with finite measure
\begin{equation}\label{DISEG}
 \P(E)\geq C \L^n(E)^{\frac{d-1}{d}}
\end{equation}
for some geometric constant $C>0$, see \cite{GN},
\cite{FGW}, and
\cite{FGuW}.
By the homogeneity properties of Lebesgue measure and $\alpha$-perimeter,
we can define the   constant
\begin{equation}\label{C_I}
     C_I=\inf\{\P(E) : \L^n(E)=1\text{ and }E\in\mathcal{S}_x,\text{ if }h\geq2\}.
\end{equation}
Only when $h\geq2$ we are adding the constraint $E\in\mathcal{S}_x$.
We have $C_I>0$ by the validity of \eqref{DISEG} for some $C>0$.
Our goal is to prove that the infimum in \eqref{C_I} is attained.

\begin{thm}
\label{thm:existence}
Let $h, k\geq 1$ and $n=h+k$. There exists an $x$-  and $y$-Schwartz
symmetric set $E\subset\Rn$ realizing
the infimum in \eqref{C_I}.

\end{thm}

\proof
Let $(E_m)_{m\in\N}$ be a minimizing sequence for the infimum
in \eqref{C_I}, with the additional assumption that
 the sets involved in the minimization
are $x$-spherically symmetric
 when $h\geq 2$.
Namely,
\begin{equation}
\label{eq:minimizing}
\L^n(E_m)=1\text{ and }\P(E_m)\leq
C_I\left(1+\frac{1}{m}\right),\quad m\in\N.
\end{equation}

By Theorems
\ref{thm:rearrangement} and \ref{thm:rear2},
we can assume that
every  set $E_m$
is $x$- and $y$-Schwartz symmetric.
We claim that the minimizing sequence
can be also assumed to be in a bounded region of $\Rn$.

Fix $m\in \N$ and let $E= E_m$.
For any $t>0$
such that \eqref{eq:density} holds 
we consider the set $E_{t-}^{x} = E\cap \{ |x|<t\}\in\mathcal S_x$.

We apply the isoperimetric inequality \eqref{DISEG}
with the constant $C_I>0$ in \eqref{C_I}
to the sets $E_{t-}^x$
and $E\setminus E_{t-}^x$,
and we use Proposition \ref{prop:derivatives}:
\begin{equation}\label{Minnie}
\begin{split}
      C_I \L^n( E_{t-}^x) ^{\frac{d-1}{d}}&\leq \P(E_{t-}^x)
    = \P(E;\{|x|<t\})+ v_E^x(t)
   \\
C_I (1-\L^n(E_{t-}^x))^{\frac{d-1}{d}}&
 \leq \P(E \setminus E_{t-}^x)
 =
\P(E; \{|x|>t\})+v^x_E(t).
\end{split}
\end{equation}
As in \eqref{paperinik2}, we let $v^x_E(t) = \mathcal{H}^{n-1}(E^x_t)$.
Adding up the two inequalities we get
\begin{equation}\label{INES}
C_I(\L^n  (E_{t-}^x)^{\frac{d-1}{d}}+(1-\L^n(E_{t-}^x))
^{\frac{d-1}{d}})\leq
\P(E)+2v_E^{x}(t).
\end{equation}
The function $g:[0,\infty)\to\R$, $g(t)=\L^n (E_{t-}^x)$ is continuous,
$(0,1)\subset g([0,\infty))\subset
[0,1]$, and it is increasing. In particular,  $g$
is differentiable almost everywhere.
For any $t>0$ such that $\P(E;\{|x|=t\})=0$, also the standard
perimeter vanishes, namely
$P(E;\{|x|=t\})=0$. With the vector field
$\varphi=(x/|x|,0)$, and for  $t<s$ satisfying
$\P(E;\{|x|=t\}) = \P(E;\{|x|=s\})=0$, we have  
\[
    \begin{split}
     \int_{E_{s-}^x\setminus E_{t-}^x}\frac{h-1}{|x|} \,dz
    & =\int_{E_{s-}^x\setminus E_{t-}^x} \mathrm{div} \varphi \,dz
     \\
& 
   = \mathcal H^{n-1}(E_s^x)-\mathcal H^{n-1}(E_t^x)+\int_{\partial^*E\cap
\{s<|x|<t\} }\langle \varphi ,\nu_E\rangle d\mathcal H^{n-1}.
    \end{split}
\]
This implies that
\[
    \lim_{s\to t} \mathcal H^{n-1}(E_s^x)=\mathcal H^{n-1}(E_t^x),
\]
with limit restricted to $s$ satisfying the above condition,
and thus
\begin{equation}
 \label{verzweigung}
      g'(t) = \lim_{s\to t}\frac{1}{s-t} \int _{t}^ s \mathcal
      H^{n-1}(E_\tau ^x)\, d\tau = \mathcal H^{n-1}(E_t^x).
\end{equation}

At this point, by \eqref{eq:minimizing},
inequality \eqref{INES} gives
\begin{equation}\label{eq:diff1}
C_I\Big(g(t)^{\frac{d-1}{d}}+(1-g(t))^{\frac{d-1}{d}}-1-\frac{1}{m}
\Big)\leq 2 g'(t).
\end{equation}
The function $\psi:[0,1]\to\R$,
$\psi(s)=s^{\frac{d-1}{d}}+(1-s)^{\frac{d-1}{d}}-1$ is concave, it
attains its maximum at $s=1/2$ with $\psi(1/2)=2^{\frac{1}{d}}-1$,
and it satisfies $\psi(s)=\psi(1-s)$, $\psi(0)=\psi(1)=0$.
By  \eqref{eq:diff1} we have
\begin{equation}\label{CIAO}
g'(t)
\geq
\frac{C_I}{2}
\Big(\psi(g(t))-\frac{1}{m}\Big)
\geq\frac{C_I}
{4}\psi(g(t))+\frac{C_I}{4}
\Big(\psi(g(t))-\frac{2}{m}\Big),
\end{equation}
for almost every $t\in\R$ and every $m\in\N$.
Provided that $m\in\N$ is such that
$2/m\leq \max\psi=2^{1/d}-1$, we show that there exist
constants $0<a_m<b_m<\infty$ such that inequality
\eqref{CIAO} implies the following:
\begin{equation}\label{eq:firstcut}
g'(t)\geq\frac{C_I}{4}\psi(g(t))\text{ for a.e.~}t\in[a_m,b_m].
\end{equation}
In fact, by continuity of $g$ and $\psi$, and by symmetry of $\psi$
with respect to the line $\{s=1/2\}$,
for $m$ large enough, there exist
$0<a_m< b_m<\infty$ such that
\[
0<g(a_m)=1-g(b_m)<\frac{1}{2}\quad
\text{and}\quad
\psi(g(a_m))=\psi(g(b_m))=\frac{2}{m}.
\]
By concavity of $\psi$
and monotonicity of $g$, it  follows that $\psi(g(t))\geq\frac{2}{m}$
for every $t\in[a_m,b_m]$, and \eqref{eq:firstcut} follows.
As $m\to\infty$ we have $g(b_{m})\to 1$, that implies
\[
 \lim_{m\to\infty} b_m = \sup\{ b>0: g(b)<1\}>0.
\]
Moreover, as $m\to \infty$ we also have $g(a_m)\to0$.
Since the set $E$ is $x$-Schwartz symmetric, there holds $g(a)>0$ for all $a>0$.
Therefore, we deduce that $a_m\to 0$.

We infer that, for $m$ large enough, we have $a_m<b_m/2$.
Integrating  inequality 
\eqref{eq:firstcut}
on the interval $[b_m/2,b_m]$, we find
\begin{equation}
\label{Pluto}
   \frac{b_m}{2}   \leq\frac{4}{C_I}\int_{b_m/2}^{b_m}\frac{g'(t)}{
\psi(g(t))}dt
\leq \frac{4}{C_I}
\int_{g(b_m/2)}^{g(b_m)}
\frac{1}{\psi(s)} ds
\leq
\frac{4}{C_I}
\int_{0}^{1}
\frac{1}{\psi(s)} ds=\ell_1.
\end{equation}

We consider the set 
$\widehat E_m=E^x_{b_m - }$. By \eqref{Pluto},
$\widehat E_m$ is contained in the
cylinder $\{|x|<2\ell_1\}$ and, 
by Proposition
\ref{prop:calibration}, it
satisfies
$
 \P(\widehat E_m)\leq \P(E_m)$.
Define the set $E_m^\dag= \delta_{\lambda_m}(\widehat E_m)$,
where
$\lambda_m\geq 1$ is chosen in such a way that $ \L^n(\widehat E_m^\dag)=1$;
namely, $\lambda_m$ is the number
\[
 \lambda_m=\Big(\frac{1}{ \L^n(\widehat E_m)}\Big)^{\frac{1}{d}},
\]
where
\begin{equation}
\begin{split}
\L^n(\widehat E_m)&= \L^n(E _m\cap\{|x|<b_m\}) = g (b_m) =1-g(a_m).
\end{split}
\end{equation}

By concavity of $\psi$, for $0<s<1/2$ the graph of $\psi$ lays above
the straight line through the origin passing through the maximum
$(1/2,\psi(1/2))$, i.e., $\psi(s)>2(2^{ {1}/{d}}-1)s$.
Therefore, since $g(a_m)<1/2$ and $\psi(g(a_m))= {2}/{m}$,
then
\[
g(a_m) \leq \dfrac{1}{m(2^{1/d}-1)},
\]
and thus
\[
  \lambda_m \leq\Big(\dfrac{1}{1-\frac{1}{m(2^{{1}/{d}}-1)}}\Big)
^{{1}/{d}}=\Big(\dfrac{m}{m-\frac{1}{2^{{1}/{d}}-1}}\Big)^{{1}/{d}}.
\]
By  homogeneity of $\alpha$-perimeter,
\begin{align*}
\P(E^\dag_m)&=\lambda_m^{d-1} \P(\hat E_m)\leq\lambda_m^{d-1}
 \P(E_m)\leq\lambda_m^{d-1}C_I\Big(1+\frac{1}{m}\Big)\\
&\leq
C_I\Big(1+\frac{1}{m}\Big)\Big(\dfrac{m}{m-\frac{1}{2^{{1}/{d}}-1}}
\Big)^{\frac{d-1}{d}}.
\end{align*}
In conclusion,   $(E^\dag_{m})_{m\in\N}$
is a minimizing sequence for $C_I$ and, for $m$ large enough,
it is contained
in the cylinder $\{|x|<\ell\}$, where $\ell=2^{1/d\,+1}\ell_1$.

Now we consider the case of the $y$-variable. We start again from
\eqref{Minnie} for the sets $E^y_{t-}$ for $t>0$.
Now the set $E$ can be assumed
to be contained in the cylinder $\{|x|<\ell\}$.
In this case, we have
\[
  v_E^y(t)=
 \int_{E_t^{y}}|x|^\alpha\,
 d\mathcal H^{n-1}
 \leq \ell^\alpha \mathcal H^{n-1} (E_t^y)
 =  \ell^\alpha g'(t).
\]
So inequality \eqref{INES} reads
\begin{equation}\label{eq:diff5}
C_I\Big(g(t)^{\frac{d-1}{d}}+(1-g(t))^{\frac{d-1}{d}}-1-\frac{1}{m}
\Big)\leq 2  \ell^\alpha g'(t).
\end{equation}
Now the argument continues exactly as in the first case.
The conclusion is that there exists a minimizing sequence
$(E_m)_{m\in\N}$ for \eqref{C_I} and there exists $\ell>0$ such that
we have:
\begin{itemize}
\item[i)] $\L^n(E_m)=1$ for all $m\in\N$;
\item[ii)] $\P(E_m) \leq C_I (1+1/m)$ for all $m\in\N$;
\item[iii)] $E_m\subset \{(x,y)\in\Rn : |x|<\ell\text{ and }|y|<\ell\}$ for all $m\in\N$;
\item[iv)] Each $E_m$ is $x$- and $y$-Schwartz symmetric.
\end{itemize}

By the compactness theorem for sets of finite $\alpha$-perimeter
(see \cite{GN} for a general statement that covers our case),
there exists   a   set $E\subset\R^{n}$
of finite $\alpha$-perimeter which is the $L^1$-limit of (a subsequence of)
the sequence  $(E_m)_{m\in\N}$.
Then we have
\[
 \L^n(E) = \lim_{m\to\infty} \L^n(E_m) = 1.
\]
 Moreover, by lower semicontinuity of
$\alpha$-perimeter
\[
   \P (E)\leq\liminf_{m\to\infty}\P(E_m)=C_I.
\]
The set $E$ is $x$- and $y$-Schwartz symmetric,
because these symmetries are preserved by the $L^1$-convergence.
This concludes the proof.
\endproof

\section{Profile  of isoperimetric sets}
\label{QUATTRO}

In Theorem \ref{thm:existence}, we proved   existence
of isoperimetric sets, in fact in the class of
 $x$-spherically symmetric sets when
$h\geq 2$.
By the characterization of the equality case
in
Theorems
\ref{thm:rearrangement} and \ref{thm:rear2},
 any isoperimetric set $E$
is $x$-Schwartz symmetric and there are
functions $c:[0,\infty)\to\R^k$
and $f:[0,\infty)\to[0,\infty)$ such that
\begin{equation}\label{Sale}
   E = \{ (x,y)\in\Rn: |y-c(|x|)|< f(|x|)\}.
\end{equation}
The function $f$ is decreasing.
We will prove in Proposition \ref{final}
that, for isoperimetric sets, the function $c$ is constant.

We start with the characterization of
an isoperimetric set $E$ with constant function $c=0$.
Let $F\subset \R^+\times \R^+$ be the generating set of $E$
\[
             E = \{(x,y)\in\Rn: (|x|,|y|) \in F\}.
\]
The set  $F$ is of the form
\begin{equation}
 \label{EFFE}
         F = \{ (r,s)\in\R^+\times \R^+: 0<s<f(r), \, r\in (0,r_0) \},
\end{equation}
where $f:(0,r_0) \to (0,\infty) $
is a decreasing function, for some $0<r_0\leq\infty$.

By the regularity theory of $\Lambda$-minimizers of perimeter,
the boundary $\partial E$ is a $C^\infty$ hypersurface
where $x\neq 0$. We do not need the general regularity theory, and we prove this fact in our case by an elementary
method that
gives also the
$C^\infty$-smoothness of
the function
$f$ in \eqref{EFFE}.

\subsection{Smoothness of $f$}\label{paperino}
We prove that the boundary $\partial F\subset\R^+\times\R^+$
is the graph of a smooth function $s=f(r)$.

We rotate  clockwise
by 45 degrees
the coordinate system $(r,s)\in\R^2$
and we call the new coordinates $(\varrho,\sigma)$; namely, we let
\[
        r = \frac{\sigma+\varrho}{\sqrt 2},\quad
        s= \frac{\sigma-\varrho}{\sqrt 2}.
\]
There exist $-\infty\leq a<0<b\leq\infty$ and a function $g:(a,b)\to\R$ such that
the boundary $\partial F\subset\R^+\times\R^+$ is a graph
$\sigma=g(\varrho)$; namely, we have
\[
       \partial F = \Big\{  \big(r(\varrho),s(\varrho)\big)=
    \Big(
\frac{g(\varrho)+\varrho}{\sqrt 2}
, \frac{g(\varrho)-\varrho}{\sqrt 2} \Big) \,: \, \varrho \in (a,b)\Big\}.
\]
Since the function $f$ is decreasing,
the function $g$ is $1$-Lipschitz continuous.

By formula \eqref{eq:N->2} and by the
standard length formula for Lipschitz graphs,
the $\alpha$-perimeter of $E$ is  
\[
     \P(E) = c_{hk} \int _a^b
\sqrt{{s'}^2 + r^{2\alpha} {r'}^2}\,  r^{h-1} s^{k-1} \, d\varrho,
\]
where $c_{hk} = hk\omega_h \omega_k $.
On the other hand, the volume of $E$ is
\[
\L^n(E)= c_{hk}
\int_a^b\left(\int_{|\varrho|}^{g(\varrho)}\left(\dfrac{\sigma+\varrho}{
\sqrt{2}}\right)^{h-1}\left(\dfrac{\sigma-\varrho}{\sqrt{2}}\right)^{k-1}
\;d\sigma \right) d\varrho.
\]
For $\varepsilon\in\R$ and $\psi\in C_c^\infty(a,b)$,
let $g_\varepsilon = g+\varepsilon
\psi$ and let $F_\varepsilon\subset\R^+\times\R^+$ be the subgraph
in $\sigma>|\varrho|$ of the function $g_\varepsilon$.
The set
  $E_\varepsilon\subset\Rn$ with generating set $F_\varepsilon$
has $\alpha$-perimeter
\[
\begin{split}
  p(\varepsilon) & =    \P(E_\varepsilon)
\\
&
= c_{hk}
\int _a^b
\sqrt{{(s'+\varepsilon\psi')}^2 + (r+\varepsilon\psi)^{2\alpha}
{(r'+\varepsilon\psi')}^2}
(r+\varepsilon \psi)^{h-1} (s+\varepsilon\psi)^{h-1} \, d\varrho,
\end{split}
\]
and volume
\[
  v(\varepsilon) = \L^n (E_\varepsilon)
 = c_{hk}
\int_a^b\left(\int_{|\varrho|}^{g(\varrho)+\varepsilon\psi(\varrho)}
 \left(\dfrac{\sigma+\varrho}{
\sqrt{2}}\right)^{h-1}\left(\dfrac{\sigma-\varrho}{\sqrt{2}}\right)^{k-1}
\;d\sigma \right) d\varrho.
\]
Since $E$ is an isoperimetric set, we have
\[
0= \left.
\dfrac{d}{d\varepsilon}\dfrac{p(\varepsilon)^d}{
v(\varepsilon)^{d-1}}\right|_{\varepsilon=0}
=\left.
\dfrac{dp^{d-1}p'v^{d-1}-p^d(d-1)v^{d-2}v'}{v^{2d-2}}\right|_{\varepsilon=0},
\]
that gives
\begin{equation}
\label{EQA}
 p'(0)-C_{hk\alpha}   v'(0)=0,\quad\text{where}\quad  C_{hk\alpha}
 =\dfrac{d-1}{d}\dfrac{\P (E)}{\L^n(E)}.
\end{equation}
After some computations, we find
\begin{equation}
\label{eq:derp}
\begin{split}
p'(0)=c_{hk}
\int_a^b\Big\{
&
\dfrac{(r^{2\alpha} r'-s') \psi'
+2\alpha r^{2\alpha-1} {r'}  ^2 \psi}
{\sqrt {{s'} ^2+ r^{2\alpha} {r'}^2 }}
+
\\
& + \sqrt{{s'} ^2
+  r^{2\alpha}{r'} ^2}
\Big[ \frac{h-1}{r} +\frac{k-1}{s}\Big]
\psi\Big\}  r^{h-1} s^{k-1}\, d\varrho,
\end{split}
\end{equation}
and
\begin{equation}
\label{eq:derv}
v'(0)=c_{hk}
\int_a^b r^{h-1}s^{k-1}\psi\,d\varrho.
\end{equation}

From \eqref{EQA}, \eqref{eq:derp}, and \eqref{eq:derv}
we deduce that $g$ is a $1$-Lipschitz function that,
via the auxiliary functions $r$ and $s$, solves
in a weak sense
the ordinary differential equation
\begin{equation}
\label{ODEg}
\begin{split}
\frac{d}{d\varrho}
\Big( r^{h-1} s^{k-1} \dfrac{r^{2\alpha} r'-s'}
{\sqrt {{s'} ^2+ r^{2\alpha} {r'}^2 }}\Big)&=
r^{h-1} s^{k-1}\Big\{
\dfrac{ 2\alpha r^{2\alpha-1} {r'}  ^2 }
{\sqrt {{s'} ^2+ r^{2\alpha} {r'}^2 }}+\\
&\quad
+ \sqrt{{s'} ^2
+  r^{2\alpha}{r'} ^2}
\Big[ \frac{h-1}{r} +\frac{k-1}{s}\Big]
- C_{hk\alpha}\Big\} .
\end{split}
\end{equation}
By an elementary argument that is omitted, if follows that
$g\in C^\infty(a,b)$.

We claim that for all $\varrho\in (a,b)$ there holds $g'(\varrho)\neq -1$.
By contradiction, assume that there exists $\bar{\varrho}\in (a,b)$
such that $g'(\bar{\varrho})=-1$, i.e., $r'(\bar{\varrho})=0$ and $s'  (\bar{\varrho})=-\sqrt{2}$.
Inserting these values into the differential equation  \eqref{ODEg}
we can compute $g''(\bar{\varrho})$ as a function of $g(\bar\varrho)$; namely, we obtain
\begin{equation}
\label{eq:secondg}
g''(\bar{\varrho})=2^{\alpha+1}
\dfrac{2(h-1)-\sqrt{2} C_{hk\alpha} r(\bar\varrho)}{r(\bar\varrho)^{2\alpha+1}}.
\end{equation}

Now there are three possibilities:
\begin{itemize}

\item[(1)] $g''(\bar\varrho)<0$.
In this case, $ g$ is  strictly concave at $\bar\varrho$
and this contradicts the fact that $E$ is  $y$-Schwartz symmetric.

\item[(2)] $g''(\bar\varrho)>0$. In this case,
   $g'$ is strictly increasing at $\bar\varrho$ and since
$g'(\bar\varrho)=-1$ this contradicts the fact the $g$ is $1$-Lipschitz,
equivalently, the fact that $E$ is $x$-Schwartz symmetric.

\item[(3)] $g''(\bar\varrho)=0$. In this case, the value of $g$ at
$\bar\varrho$ is, by \eqref{eq:secondg},
\begin{equation}
 \label{picco}
   g(\bar\varrho)=-\bar\varrho+\dfrac{\sqrt{2}(h-1)}{C_{hk\alpha}}.
\end{equation}
The function
$  \widehat g(\varrho)=-\varrho
+\dfrac{\sqrt{2}(h-1)}{C_{hk\alpha}}$, $\varrho\in\R$,
is the unique solution to the ordinary differential equation
\eqref{ODEg} with initial conditions $g(\bar\varrho)$ given by \eqref{picco}
and
$g'(\bar\varrho)=-1$. It follows that $g=\widehat g$ and
this contradicts the boundedness of the isoperimetric set; namely,
the fact that isoperimetric sets have finite measure.
\end{itemize}

This proves that $g'(\varrho)\neq -1$ for all $\varrho\in(a,b)$.

\subsection{Differential equations for the profile function}\label{paperino2}

By the discussion in the previous section, the function
$f$ appearing in the definition of the set $F$ in \eqref{EFFE}
is in $C^\infty(0,r_0)$. The function $f$ is decreasing, $f'\leq 0$.
By formula \eqref{eq:N->2}, the perimeter
of the set $E$ with generating set $F$ is
\begin{equation}
\label{eq:perflip}
\begin{split}
\P (E)= c_{hk}
\int_{0}^{r_0}  \sqrt{f'(r)^2+r^{2\alpha}}\;r^{
h-1}f(r)^{k-1}dr,
\end{split}
\end{equation}
and
the volume of $E$ is
\begin{equation}
\label{eq:volumflip}
\L^n(E)= \frac{c_{hk}}{k}  \int_0^{r_0} r^{h-1}f(r)^k dr.
\end{equation}
As in the previous subsection, for   $\psi\in C_c^\infty(0,r_0)$
and $\varepsilon\in\R$, we consider the
perturbation $f+\varepsilon\psi$
and we define the set
\[
 E_\varepsilon = \big\{(x,y) \in\Rn :  |y| < f(|x|)+\varepsilon \psi(|x|)\big\}.
\]
Then we have
\begin{equation*}
\begin{split}
 p(\varepsilon) & =\P(E_\varepsilon)=c_{hk}
 \int_0^{r_0}  \sqrt{
(f'+\varepsilon\psi')^2+r^{2\alpha}}\,
(f+\varepsilon\psi)^{k-1}r^{h-1}
dr,
\\
v(\varepsilon)&
=\L^n(E_\varepsilon)=
\frac{c_{hk}}{k}\int_0^{r_0}
(f+\varepsilon\psi)^{k}r^{h-1}dr,
\end{split}
\end{equation*}
and from these formulas we compute
the first derivatives at $\varepsilon=0$:

\begin{equation*}
\begin{split}
   p'(0) & =c_{hk}
 \int_0^{r_0}\Big[\dfrac{f^{k-1}f'}{\sqrt{
{f'}^2+r^{2\alpha}}}\psi'+(k-1)f^{k-2}\sqrt{f'^2+r^{2\alpha}}
\psi\Big]\, r^{h-1}\, dr,
\\
v'(0) & = c_{hk}  \int_0^{r_0} f^{k-1}\psi \, r^{h-1}\, dr.
\end{split}
\end{equation*}
The minimality equation \eqref{EQA} reads
\begin{equation}
\label{eq:necessary1}
\int_0^{r_0}
\Big(\dfrac{f'f^{k-1}}{\sqrt{f'^2+r^{2\alpha}}}
\psi'+\Big[(k-1)f^{k-2}\sqrt{f'^2+r^{2\alpha}}-C_{hk\alpha}
 f^{k-1}\Big]
\psi\Big)\, r^{h-1}\, dr=0.
\end{equation}
Integrating by parts
the term with $\psi'$
and using the fact that $\psi$ is arbitrary,
we deduce that $f$ solves the following second order ordinary differential equation:

\begin{equation}
\label{eq:diffnecessary}
-\dfrac{d}{dr} \left(r^{h-1}\dfrac{f'f^{k-1}}{\sqrt{f'^2+r^{
2\alpha}}}\right)+r^{h-1}\left[(k-1)\sqrt{f'^2+r^{2\alpha}}\;f^{k-2
}-C_{hk\alpha} f^{k-1}\right]=0.
\end{equation}
The normal form of this differential equation is

\begin{equation}
\label{eq:odef}
f''=\dfrac{\alpha
f'}{r}+({f'}^2+r^{2\alpha})\left(\dfrac{k-1}{f}-(h-1)\dfrac{f'}{r^{
2\alpha+1}}\right)-C_{hk\alpha}
 \dfrac{({f'}^2+r^{2\alpha})^{\frac{3}{2}}}{r^{2\alpha}
},
\end{equation}
and it  can be rearranged in the following ways:
\begin{equation}\label{pongo}
\begin{split}
\dfrac{\partial}{\partial r}\Big(
\dfrac{f'}{r^\alpha}\Big)
&
=   ({f'}^2+r^{2\alpha})
\Big(\dfrac{k-1}{f r^\alpha }-(h-1)\dfrac{f'}{
r^{3\alpha+1}}\Big)-C_{hk\alpha}
\dfrac{(f'^2+r^{2\alpha})^{\frac{3}{2}}}
{r^{3\alpha}}
\\
&
=r^{\alpha}\Big(\Big(
\dfrac{f'}{r^\alpha}\Big)^2+1\Big)
\Big(\dfrac
{k-1}{f}-\dfrac{(h-1)}{r^{\alpha+1}}\dfrac{f'}{r^{\alpha}}
\Big)-C_{hk\alpha}\Big(\Big(\dfrac{f'}
{r^{\alpha}}\Big)^2+1\Big)^{\frac{3}{2}}.
\end{split}
\end{equation}

With the substitution
\begin{equation}\label{ZETA}
z=\sin\arctan\Big(\frac{f'}{r^\alpha}\Big)=\frac{f'}{\sqrt{r^{2\alpha}+f'^{2}}},
\end{equation}
equation \eqref{pongo} transforms into the equation
\begin{equation}\label{pop}
(r^{h-1} z)' =   r^{\alpha+h-1}  \dfrac{k-1}{f}
   \sqrt{1-z^2}-C_{hk\alpha} r^{h-1}.
\end{equation}
We integrate this equation on the interval $(0,r)$. When $h>1$ we use
the fact that $r^{h-1} z =0 $ at $r=0$. When $h=1$ we use the fact that $z$ has
a finite limit as $r\to0^+$. In both cases, we  deduce that there exists
a constant $D\in\R$ such that
\begin{equation}\label{EZETA}
    z(r)=r^{1-h}\int_0^rs^{\alpha+h-1}\frac{k-1}{f}\sqrt{1-z^2}\;ds-\frac{C_{
hk\alpha}}{h}r+Dr^{1-h}.
\end{equation}
Inserting \eqref{ZETA} into \eqref{EZETA}, we get
\begin{equation}\label{zorro}
   \frac{f'}{\sqrt{r^{2\alpha}+f'^{2}}} = r^{1-h}  \int _0 ^r s^{2\alpha+h-1}  \dfrac{k-1}{f\sqrt{s^{2\alpha}+f'^{2}}}
    \, ds -\frac{C_{hk\alpha} }{h}r +D r^{1-h}.
\end{equation}
If $h\geq2$, from \eqref{zorro} we deduce that $D=0$.
In fact,  the left-hand side of \eqref{zorro} is bounded
as $r\to0^+$, while the right-hand side diverges to $\pm\infty$
according to the sign of $D\neq 0$.
In the next section, we prove that $D=0$ also when $h=1$,
provided that  $f$ is the profile of an  
isoperimetric set.

\medskip

\begin{rem}[Computation of the solution when  $k=1$] \label{sol1} When $k=1$ and $D=0$,
equation 
\eqref{zorro} reads
\[
   \frac{f'}{\sqrt{r^{2\alpha}+f'^{2}}} = -\frac{C_{hk\alpha} }{h}r.
\]
and this is equivalent to
\begin{equation}\label{pimpa}
f'(r) = -\frac{ C_{hk\alpha} r^{\alpha+1}}{\sqrt{h^2 - C_{hk\alpha}^2   r ^2}},\quad r\in
[0,r_0).
\end{equation}
Without loss of generality we can assume that $r_0=1$
and this holds if and only if  ${C_{hk\alpha}}={h}$.
Integrating \eqref{pimpa} with $f(1)=0$ we
obtain the solution
\[
 f(r) = \int_r^1 \frac{s^{\alpha+1} }{\sqrt{1-s^2}}\, ds = \int_{\arcsin
r}^{\pi/2}\sin^{\alpha+1}(s) \,
ds.
\]
This is the profile function for the isoperimetric
set when $k=1$ in \eqref{EXPO}.
\end{rem}

\subsection{Proof that $D=0$ in \eqref{zorro}.}\label{paperino3}

We prove that $D=0$ in the case $h=1$.
We assume by contradiction that $D\neq0$.
For a small parameter $s>0$, let
 $f_s:[0,r_0)\to \R^+$ be the function
\[
 f_s(r)=\left\{\begin{array}{ll}f(s)&\text{ for
} 0 < r \leq s
\\
f(r)
&
\text{ for }r>s,
\end{array}
\right.
\]
and define the set
\[
   E_s = \big\{ (x,y)\in \Rn :  |y|< f_s(|x|)\big\}.
\]
Recall that the isoperimetric ratio is $
\I_\alpha(E)={\P(E)^d}/{\mathcal L^n(E)^{d-1}}$.
We claim that for $s>0$ small, the difference of isoperimetric ratios
\begin{equation}
\label{CLAM}
\begin{split}
\I_\alpha(E_s) - \I_\alpha(E) & = 
\dfrac{\P(E_s)^d}{\L^n(E_s)^{d-1}}-\dfrac{\P(E)^d}{\L^n(E)^{d-1}}
\\
&
=
\dfrac{\P(E_s)^d\L^n(E)^{d-1}-\P(E)^d\L^n(E_s)^{d-1}}{\L^n(E_s)^{d-1}
\L^n(E)^{d-1}}
\end{split}
\end{equation}
is strictly negative.

The  $\alpha$-perimeter
of $E_s$ is
\begin{align*}
\P(E_s) & =
 c_{hk}  \int_0^\infty\sqrt{{f'_s}^2+r^{2\alpha
}}f_s^{k-1} \, r^{h-1}\,dr
\\
&=c_{hk} \Big[
 f(s)^{k-1} \int_0^s r^{\alpha+h-1} \;dr
+\int_s^\infty\sqrt{f'^2+r^{2\alpha}}f^{k-1}\, r^{h-1}
\,dr\Big]
\\
&=\P (E)+c_{hk}
\int_0^s\Big[r^\alpha
f(s)^{k-1}-\sqrt{f'^2+r^{2\alpha}}f^{k-1}\Big]r^{h-1}\,dr,
\end{align*}
and its   volume is
\begin{align*}
\L^n(E_s)  &= \frac{c_{hk}}{k}
 \int_0^\infty f_s^k r^{h-1}\,dr
=\frac{c_{hk}}{k}
 \Big(\int_0^s f(s)^k r^{h-1}\,dr
+ \int_s^\infty f(r)^k\, r^{h-1}\,dr\Big)
\\
&=\L^n(E)+\frac{c_{hk}}{k}\int_0^s
\big(f(s)^k-f(r)^k\big)\, r ^{h-1}\,dr,
\end{align*}
so, by elementary Taylor approximations, we find
\[
\begin{split}
\L^n&(E)  ^{d-1} \P(E_s)^d
 =
\\
&
= \L^n(E)^{d-1}
 \Big\{\P (E)+c_{hk}
 \int_0^s
 \Big[ r^\alpha f(s)^{k-1}-\sqrt{f'^2+r^{2\alpha}}
 f^{k-1}\Big]\, r ^{h-1}\;dr\Big\}^d
\\
&
= \L^n (E)^{d-1} \bigg\{
\P(E)^d+ d c_{hk} \P (E)^{d-1}
\bigg(\int_0^s\big[r^\alpha f(s)^{k-1}
-\sqrt{f'^2+r^{
2\alpha}}f^{k-1}\big]\, r^{h-1}\,d r \bigg)
\\
&
\qquad\qquad\qquad\qquad
\qquad\qquad
+R_1 (s)
\bigg\},
\end{split}
\]
where $R_1(s)$ is a higher order infinitesimal as $s\to0$, and
\[
\begin{split}
\P(E)^{d}&
\L^n(E_s)^{d-1}  =
\P(E)^d \Big\{\L^n(E)+ \frac{c_{hk}}{k}
\int_0^s
\big(f(s)^k-f(r)^k\big)\, r^{h-1}\,dr\Big\}^{d-1}
\\
&
=
\P(E)^d
 \Big\{\L^n(  E)^{d-1}  + \frac{c_{hk} (d-1)}{k}  \L^n(E)^{d-2}
\int_0^s \big(f(s)^k-f(r)^k\big)\, r^{h-1}\, dr+R_2(s)  \Big\},
\end{split}
\]
where $R_2(s)$ is a higher order infinitesimal as $s\to0$.
The difference is thus
\[
 \begin{split}
  \Delta(s) &= 
P(E_s)^d\L^n(E)^{d-1}-\P(E)^d\L^n(E_s)^{d-1}
\\
 &
 =
c_{hk}\P(E)^d
 \L^n(E)^{d-1}
\Big\{
d  \frac{A(s)}{\P(E)}
 -(d-1)
 \frac{B(s)}{k \L^n(E)}\Big\},
 \end{split}
\]
where we let
\[
\begin{split}
A(s) &
 =  \int_0^s\big[r^\alpha f(s)^{k-1} -\sqrt{f'^2+r^{
2\alpha}}f^{k-1}\big]\, r^{h-1}\,d r +R_1(s)
\\
B(s) & = \int_0^s \big(f(s)^k-f(r)^k\big)\, r^{h-1}\, dr+R_2(s).
\end{split}
\]

Now we let $h=1$ and we observe that the differential equation \eqref{EZETA} or
its equivalent
version \eqref{zorro} imply that 
\[
 \lim_{r\to0^+} \frac{f'(r)}{r^\alpha} = D.
\]
So for $D\neq 0$ and, in fact, for  $D<0$ (because $f$ is decreasing) we have
\[
 \lim_{s\to0^+} 
\frac{A(s)}{s^{\alpha+h}} 
=
f(0)^{k-1} 
\frac{ 1-\sqrt{D^2+1} }{\alpha+h}<0,
\]
and
\[
  \lim_{s\to0^+} 
\frac{B(s)}{s^{\alpha+h}}=0.
\]
It follows that for $s>0$ small there holds
\[
 \frac{ \Delta(s)}{ s^{h+\alpha}} = 
f(0)^{k-1} 
\frac{ 1-\sqrt{D^2+1} }{\alpha+h}d c_{hk}\P(E)^{d-1}
 \L^n(E)^{d-1} 
 +o(1) <0.
\]
Then $E$ is not an isoperimetric set. This proves that $D=0$.

\subsection{Initial  and final conditions for the profile function}
In this section, we study the behavior of $f$ at $0$ and $r_0$.

\begin{prop}\label{propo2}
The profile function $f$ of an $x$- and $y$-Schwartz symmetric
isoperimetric set $E\subset\Rn$
satisfies  $f \in C^\infty(0,r_0)\cap C([0,r_0])$  for some $0<r_0<
\infty$,   $f'\leq 0$,  $f(r_0) = 0$, it  solves   the
differential equation
\eqref{zorro} with $D=0$, and
\[
\lim_{r\to r_0^-} f'(r) = -\infty \quad \textrm{and}
\quad
 \lim_{r\to 0^+} \frac{f'(r)}{ r^{\alpha+1}} = -\frac{C_{hk\alpha}
}{h}.
\]
\end{prop}

\proof 
By Remark \ref{sol1}, it is sufficient to prove that $r_0<\infty$ when $k>1$.
Assume by contradiction that $r_0=\infty$. In this case, it must be
\begin{equation}\label{limme}
 \lim_{r\to\infty} f(r) = 0,
\end{equation}
otherwise the set $E$ with profile $f$ would have infinite volume.

For $\varepsilon >0$ and $M>0$, let us consider
the set
\[
 K_M =\big\{ r\geq M: f'(r)\geq -\varepsilon\big\}.
\]
Recall that in our case we have $f'\leq 0$.
The set $K_M$ is closed and nonempty for any $M$. If $K_M=\emptyset$ for some
$M$,
then this would contradict \eqref{limme}.

Let $\bar r\in K_M$.
From \eqref{eq:odef} we have
\begin{equation}\label{secco}
\begin{split}
  f''(\bar r) &
 = - \frac{\alpha\varepsilon}{\bar r}+\bar r^{2\alpha}
\frac{k-1}{f(\bar r)} - C_{hk\alpha} \frac{(\varepsilon^2+ r\bar
r^{2\alpha})^{3/2}}{\bar r^{2\alpha}}
\\
&
 \geq \frac 12 M^{2\alpha} \frac{k-1}{f(M)}>0,
\end{split}
\end{equation}
provided that $M$ is large enough.
We deduce that there exists $\delta>0$ such that
$f'(r)\geq -\varepsilon$ for all $r\in [\bar r,\bar r+\delta)$.
This proves that $K_M$ is open to the right. It follows that it must be $K_M =
[M,\infty)$. This proves that 
\[
 \lim_{r\to\infty} f'(r) = 0,
\]
and this in turn contradicts \eqref{secco}.

Now we have $r_0<\infty$ and we also have
\[
 L = \lim_{r\to r_0^-} f(r) = 0.
\]
If it were $L>0$, then
the isoperimetric set would have a ``vertical part''.
We would get a contradiction by the argument at  point (3) 
at the end of Section  
\ref{paperino}.

We claim that 
\[
 \lim_{r\to r_0^-} f'(r) = -\infty.
\]
For $M>0$ and $0<s<r_0$, consider the set
\[
 K_s = \big\{ s\leq r< r_0: f'(r)\geq -M\big\}.
\]
By contradiction assume that there exists $M>0$ such that
$K_s\neq \emptyset$ for all $0<s<r_0$. If $\bar r \in K_s$, we have
as above $f''(\bar r) \geq \frac{1}{2}(k-1) s^{2\alpha} /  f(s)>0$. We deduce
that there exists $s<r_0$ such that $0\geq  f'(r) \geq -M$ for all
$r\in[s,r_0)$. From \eqref{eq:odef}, we deduce that there exists a constant
$C>0$
such that  
\[
 f''(r) \geq \frac{C}{f(r)}.
\]
Multiplying by $f'\leq 0$ and integrating the resulting inequality we find
\[
  f'(r) ^2 \leq 2C \log | f(r)| +C_0,
\]
for some constant $C_0\in\R$.
This is a contradiction because $\displaystyle \lim_{r\to r_0^-} \log |f(r)|
=-\infty$.

By Section \ref{paperino2}, we have   $D=0$ in \eqref{zorro}.
In this case,  by \eqref{EZETA} we can compute the limit
\[
 \begin{split}
\lim_{r\to0^+} \frac{f'(r)}{ r^{\alpha +1}} =
\lim_{r\to0^+}  -\frac{C_{hk\alpha} }{h}+
r^{-h}  \int _0 ^r s^{\alpha+h-1}  \dfrac{k-1}{f}
   \sqrt{1-z^2} \, ds
=-\frac{C_{hk\alpha} }{h}.
\end{split}
\]
This ends the proof.
\endproof

\begin{rem} The Cauchy Problem for the differential equation
\eqref{eq:odef}, with the initial conditions
$f(0)=1$ and $f'(0)=0$ 
has a unique decreasing solution on some interval $[0,\delta]$, with $\delta>0$,
in the class of functions $f\in C^1([0,\delta])\cap C^\infty((0,\delta])$  
such that 
\[
\lim_{r\to0^+}\dfrac{f'(r)}{r^{\alpha+1}}=-\dfrac{C_{hk\alpha}}{h}.
\]
This can be proved using the Banach fixed point Theorem with the norm
\[
\|f\|=\max_{r\in[0,\delta]}|f(r)|+\max_{r\in[0,\delta]}\dfrac{|f'(r)|}{r^{\alpha+1}}.
\]

From Theorem \ref{thm:existence} and Proposition \ref{propo2}, there exists
a value of the constant $C_{hk\alpha}>0$ such that the maximal decreasing solution
of the Cauchy Problem has  a maximal interval $[0,r_0]$ such that $f(r_0)=0$.

\end{rem}

\subsection{Isoperimetric sets are $y$-Schwartz symmetric}
To conclude the proof of Theorems \ref{THM1}--\ref{THM3}
we are left to show
that for an isoperimetric set $E$ of the type \eqref{Sale},
the function $c$ of the centers is constant.

\begin{prop}
\label{final}
Let $h,k\geq 1$ and $n=h+k$. Let
$E\subset\Rn$ be a set of the form 
\[
  E = \{ (x,y) \in \Rn: |y-c(|x|)|< f(|x|)\}
\]
for measurable  
functions $c:[0,\infty) \to\R^k$ and
$f:[0,\infty)\to[0,\infty]$.
If $E$ is an isoperimetric set for the problem \eqref{C_I} then the function $c$ is constant.

\end{prop}

\proof
If $E$ is isoperimetric, then  also its $y$-Schwartz
rearrangement
$E^* = \{ (x,y) \in \Rn: |y|< f(|x|)\}$ is an isoperimetric set,
see Theorems \ref{thm:rearrangement} and \ref{thm:rear2}. Then, by
Proposition \ref{propo2}, we
have $  f \in C^\infty(0,r_0)\cap C([0,r_0])$ with $f(r_0)=0$ and $f'\leq 0$. In particular, $f\in \mathrm{Lip}_{\mathrm{loc}}(0,r_0)$.
We claim that  $c\in \mathrm{Lip}_{\mathrm{loc}}(0,r_0)$.

Since $E$ is $x$-Schwartz symmetric, for any $0<r_1<r_2<r_0$ we have the inclusion
\[
 \{ y\in \R^ k : |y- c(r_2)| \leq f(r_2)\}\subset  
 \{ y\in \R^ k : |y- c(r_1)| \leq f(r_1)\}.
\]
Assume $c(r_2)\neq c(r_1)$ and let $\vartheta = c(r_2) -c(r_1)/|c(r_2)
-c(r_1)|$. Then we have
\[
    c(r_2) +\vartheta f(r_2) \in 
 \{ y\in \R^ k : |y- c(r_1)| \leq f(r_1)\},
\]
and therefore
\[
|c(r_2) - c(r_1)| + f(r_2) =        | c(r_2) +\vartheta f(r_2) - c(r_1)| \leq
f(r_1).
\]
This implies that $c$ is locally Lipschitz on $(0,r_0)$.

Let $F\subset \R^+\times \R^k$ be the generating set of $E$:
\[
    E = \{(x,y) \in \R^n : (|x|, y) \in F\}.
\]
By the discussion above, the set $E$ and thus also the set $F$ have locally
Lipschitz boundary away from a negligible set.
By the representation formula
\eqref{eq:N->k+1}, we have
\[
 \P(E) = Q(F) = h\omega_h \int_{\partial F} \sqrt{N_r^2 + r^{2\alpha}
|N_y|^2}\, r^{h-1} \, d\mathcal H^{k},
\]
where $(N_r, N_y)\in \R^{1+k} $ is the unit normal to $\partial F$ in
$\R^+\times\R^k$, that is defined $\mathcal H^{k}$ almost everywhere on the
boundary.
By the coarea formula (see \cite{BZ}) we also have 
\[
 Q(F) =h\omega_h  \int _0^\infty r^{h-1} \int _{\partial F_r }\frac{ \sqrt{N_r^2
+ r^{2\alpha}
|N_y|^2}}{\sqrt{1-N_r^2}}\, d\mathcal H^{k-1}\, dr ,
\]
where $\partial F _r = \partial \{ y\in \R^k: (r,y) \in F\} = \{y\in\R^k :
|y-c(r)| = f(r)\}$.

A defining equation for $\partial F$ is $|y-c(r)|^2 - f(r)^2 =0$.
From this equation, we find
\[
\begin{split}
 N_r  & = -\frac{\langle y-c,c'\rangle +ff'}{\sqrt{(\langle
y-c,c'\rangle +ff')^2+ |y-c|^2}},
\\
 N_y   & = \frac{y-c}{\sqrt{(\langle y-c,c'\rangle
+ff')^2+ |y-c|^2}},
\end{split}
\]
and thus, by translation and scaling in the inner integral,
\[
\begin{split}
 Q(F) &
   = h \omega_ h \int _0^\infty r^{h-1} \int _{|y-c(r)|=f(r)}
\sqrt{\Big\{\frac{\langle y-c(r),c'(r)\rangle}{f(r)}+ f'(r)\Big\}^2
+r^{2\alpha} 
}\, d\mathcal H^{k-1}(y)\, dr
\\
& 
   = h \omega_ h \int _0^\infty r^{h-1} f(r)^{k-1} \int _{|y|=1}
\sqrt{\big\{ \langle y,c'(r)\rangle+ f'(r)\big\}^2 +r^{2\alpha} 
} \, d\mathcal H^{k-1}(y)\, dr. 
\end{split}
\]
For any $r>0$, the function $\Phi: \R^h\to \R^+$
\[
 \Phi(z) = \int _{|y|=1}
\sqrt{\big(\langle y,z\rangle+ f'(r)\big)^2 +r^{2\alpha} 
} \, d\mathcal H^{k-1}(y)
\]
is strictly convex. This follows from the strict convexity of
$t\mapsto\sqrt{r^{2\alpha}+t^2}$.
The function $\Phi$ is also  radially symmetric because the integral
is invariant under orthogonal transformations. It follows that $\Phi$ 
attains the minimum at the point $z=0$ and that this minimum point is unique.

Denoting by $F^*$   the generating set  of $E^*$, we deduce that
if $c'$ is not $0$ a.e., then we have the strict inequality
$\P(E^*) = Q(F^*)<Q(F)=\P(E)$, and $E$ is not isoperimetric.
Hence, $c$ is constant and
this concludes the proof.

\endproof

\end{document}